\documentclass[a4paper]{article}

\usepackage[a4paper,
            bindingoffset=0.2in,
            left=1in,
            right=1in,
            top=1in,
            bottom=1in,
            footskip=.25in]{geometry}

\usepackage{graphicx}
\usepackage{float}
\usepackage{subfigure}
\usepackage{amsfonts, amsmath, amsthm, amssymb}
\usepackage{orcidlink}
\usepackage{authblk}
\usepackage{tikz}
\usepackage{pgfplots}
\pgfplotsset{compat=newest}
\usepackage[maxbibnames=9]{biblatex} 
\usepackage{hyperref}
\usepackage{comment}
\usepackage{mathtools}
\usepackage{todonotes}

\usepackage{algorithm}
\usepackage{algpseudocode}

\newtheorem{theorem}{Theorem}[section]
\newtheorem{lemma}[theorem]{Lemma}
\newtheorem{definition}[theorem]{Definition}
\newtheorem{remark}[theorem]{Remark}
\newtheorem{proposition}[theorem]{Proposition}
\newtheorem{corollary}[theorem]{Corollary}
\newtheorem*{nat}{Notation and terminology}

\def\R{{\mathbb R}}

\def\P{{\mathbb P}}

\newcommand{\bbP}{\mathbb{P}}

\newcommand{\de}{\,\mathrm{d}}

\newcommand{\spa}{\mathrm{span} }

	{\left\lbrace \begin{array}{@{} l @{} }}%
	{ \end{array}\right.}

\ifpdf
  \DeclareGraphicsExtensions{.eps,.pdf,.png,.jpg}
\else
  \DeclareGraphicsExtensions{.eps}
\fi

\usepackage{xcolor}
\definecolor{turchese}{RGB}{35, 174, 163}

\newcommand{\lbb}[1]{\textcolor{black}{#1}}

\title{Polynomial approximation from diffused data: unisolvence and stability} 

\author[a,b]{Ludovico Bruni Bruno\footnote{email \href{mailto:ludovico.brunibruno@polito.it}{ludovico.brunibruno@polito.it}} \orcidlink{0000-0002-5246-8049}}

\affil[a]{\footnotesize Dipartimento di Scienze Matematiche \textquotedblleft Giuseppe Luigi Lagrange\textquotedblright, Politecnico di Torino, corso Duca degli Abruzzi, 24, Torino, 10129, Italia}

\author[b]{Stefano De Marchi\footnote{email \href{mailto:stefano.demarchi@unipd.it}{stefano.demarchi@unipd.it}}
\orcidlink{0000-0002-2832-8476}}

\affil[b]{\footnotesize Dipartimento di Matematica \textquotedblleft Tullio Levi-Civita\textquotedblright, Università di Padova, via Trieste, 63, Padova, 35131, Italia}
            
\author[c]{Giacomo Elefante\footnote{email \href{mailto:giacomo.elefante@usi.it}{giacomo.elefante@usi.it}} \orcidlink{0000-0001-5576-6802}}

\affil[c]{\footnotesize Istituto Dalle Molle di studi sull'intelligenza artificiale, Università della Svizzera italiana, via La Santa, 1, Lugano, 6962, Svizzera}
\date{}

\begin{document}

\maketitle

\begin{abstract}
In this work, we address the problem of polynomial interpolation of non-pointwise data. More specifically, we assume that our input information comes from measurements obtained on diffuse compact domains. Although the nodal and the diffused problems are related by the mean value theorem, such an approach does not provide any concrete insights in terms of well-posedness and stability. We hence develop a different framework in which {\it unisolvence} can be again recovered from nodal results, for which a wide literature is available. 
To analyze the stability of the so-obtained diffused interpolation procedure, we characterize the norm of the interpolation operator in terms of a Lebesgue constant-like quantity. After analyzing some of its features, such as invariance properties and sensitivity to support overlapping, we numerically verify the theoretical findings.
\end{abstract}

\section{Introduction}

Interpolation deals with the reconstruction and approximation of a function $ f(x) $ given some samplings $ f(x_1), \ldots, f(x_N) $ on a unisolvent set $ x_1, \ldots, x_N $ of nodes in some domain $ \Omega \subset \R^n $. This procedure assumes that $ f $ admits a point-wise value at such nodes and, in a problem-oriented perspective, that we are also able to capture the corresponding data. In some situations, for instance in climate studies and geophysics \cite{DynWahba}, physical \cite{WangJCP} or statistical \cite{Schoenberg} problems, optimal transport \cite{OptimalTransport}, mimetic methods \cite{Hiemstra} or sensor placement \cite{SensorPlacement}, such data is not available or may even be meaningless, e.g. if we are trying to interpolate a differential $n$-form \cite{RapettiAubry} by means of weights \cite{Rapetti07}. This latter situation arises in the construction of finite element spaces \cite{ABR23,BBZ22} \lbb{and appears in virtual elements related settings \cite{ASV20}}. In such cases, \emph{diffused} data can be in fact given. This forces to replace of the nodal interpolation procedure with \emph{integral interpolation} \cite{IntegralInterpolation}. As a result, the usual interpolating conditions that allow to project a function $ f $ onto a (total) degree $ d $ polynomial $ p(x) \in \P_d (\Omega) $ are replaced by
\begin{equation} \label{eq:areamatching}
\int_{K_i} f(x) \de x = \int_{K_i} p(x) \de x , \qquad i = 1, \ldots, N := \dim \P_d (\Omega), 
\end{equation}
or their normalized (by the Lebesgue measure of $ K_i$, henceforth denoted by $ \mu (K_i)  $) variants \cite{BeiraoMimetic}, being $ \{ K_i \}_{i=1}^N $ a collection of compact sets supported in $ \Omega $. 
\begin{remark} \label{rmk:assumptionscompacts}
    Any compact set $ K $ in the present paper is understood to be the closure of open and connected set $ U $ of finite and non-vanishing Lebesgue measure: $0 < \mu(K) = \mu(U) < +\infty $. 
\end{remark}
Conditions in \eqref{eq:areamatching} are called \emph{area-matching} \cite{Schoenberg} and the corresponding technique is known as \emph{histopolation} \cite{Robidoux}, as (at least when the dimension of the ambient space is $ n = 1 $) they correspond to asking that the area spanned by some histograms associated with $ f(x) $ and $ p (x) $ coincide. When supports $K_i$ are tensor products of segments, several spline-based techniques have been developed, see e.g. \cite{Duchon}, but for less regular domains, histopolation still has many obscure sides.

Now, suppose that indeed the information at the left hand side of \eqref{eq:areamatching} represent some data extracted from a sensor. For such an object, we may decide both its shape and its placement in the domain $ \Omega $. Of course, we want to place these $ N $ objects in $ \Omega $ in such a way that not only is the problem \eqref{eq:areamatching} well posed (that is, for each $ f(x) $ there exists a unique polynomial $ p(x) $ of total degree at most $ d $ satisfying such a collection of degrees of freedom), but also that the distortion of the reconstructed data is as small as possible, for instance avoiding Runge phenomena \cite{BruniRunge} or, more generally, limiting uncontrollable propagation of instability in an appropriate sense.

The well-posedness of the interpolator $ \Pi: C^0 (\Omega) \to \P_d (\Omega) $ associated with \eqref{eq:areamatching} is related to the \emph{unisolvence} of the supports $K_i$. 

\begin{definition}
    A set of compacts $ \{ K_1, \ldots, K_N \} $ is termed unisolvent for $ \P_d (\Omega) $ if $ \int_{K_i} p(x) \de x = 0 $ for $ i = 1, \ldots, N $ implies that $ p (x) = 0 $.
\end{definition}

The characterization of unisolvent sets in the present framework is far from being trivial. We face this problem in Section \ref{sect:unisolvence}, where we devise two techniques for the identification of unisolvent sets. 
It is worth mentioning that such techniques hinge on the interplay with nodal unisolvence, for which a rich literature is available. Without the presumption of being exhaustive enough, we recall the papers \cite{InterpolationOnFlats,Bos91,Calvi,Guenter} and some of their relevant applications \cite{Chkifa,CY77,Marchildon}. An approach that well matches our construction is curve-based one, considered for instance in \cite{PaduaCurve,ErbLissajous,Floater}.

Once unisolvence is acquired, it is possible to select a convenient basis $ \{ \ell_{K_1}, \ldots, \ell_{K_N} \} $ for $ \P_d (\Omega) $. 

\begin{definition}
    The Lagrange basis for $ \P_d (\Omega) $ is the unique basis satisfying the duality relations
\begin{equation} \label{eq:Lagrangeduality}
    \int_{K_i} \ell_{K_j}(x) \de x = \delta_{i,j} \qquad i,j = 1, \ldots, N .
\end{equation}
\end{definition}

The Lagrange basis allows us to represent the solution of the problem \eqref{eq:areamatching} in a convenient form
\begin{equation} \label{eq:Lagrangerepresentation}
    p(x) = \Pi f(x) = \sum_{i=1}^N \left( \int_{K_i} f(x) \de x \right) \ell_{K_i} (x) ,
\end{equation}
which in turn also shows that the diffused interpolator is also a projector, see Proposition \ref{prop:projector}.

The representation \eqref{eq:Lagrangerepresentation} separates the contribution of the specific problem (carried by $ \int_{K_i} f $), which is here used as a coefficient on the polynomial expansion, and the impact of the choice on the supports $ K_i $, determining the basis functions $ \ell_{K_i}$. To have a measure of the quality of our supports, we shall hence choose compacts that yield Lagrange functions with small sup-norm \cite{Davis}; recall that
$$ \Vert f \Vert_\infty := \sup_{x \in \Omega} | f(x) | .$$
That this norm might be meaningful also in our approach is suggested by the mean value theorem applied to the set of conditions in Eq. \eqref{eq:areamatching}, as it transforms the diffused interpolation problem into a nodal one $ f(\xi_i) = p(\xi_i) $ for $ i = 1, \ldots, N $, on some \emph{unknown} points $ \xi_i \in K_i $ for $ i = 1, \ldots, N $.
To relate the uniform norm with the diffused supports $ K_i $, we may again invoke the mean value theorem to see that
\begin{equation} \label{eq:equivalencenorm}
    \Vert f \Vert_\infty := \sup_{x \in \Omega} | f(x) | = \sup_{K \subset \Omega} \frac{1}{\mu(K)} \left| \int_K f(x) \de x \right| .
\end{equation}
After recalling some classical inequalities that enlighten the role of the operator norm $ \Vert \Pi \Vert_{\mathrm{op}} $ of $ \Pi $ in uniform approximation, in Section \ref{sect:uniform} we exploit \eqref{eq:equivalencenorm} to relate $ \Vert \Pi \Vert_{\mathrm{op}} $ with the underlying diffused supports. This provides a Lebesgue constant \cite{ARR20}, thus a measurement of the numerical conditioning, for the diffused problem.

Once the features of this Lebesgue constant are established, it is necessary to detect collections of supports that keep such a quantity under control. As in the nodal framework, this Lebesgue constant turns out to be extremely sensitive to the choice of the supports. The identification of effective sets of diffused supports is also complicated by the large number of parameters needed to represent a compact (in place of a point), which generally leads to an undetermined problem. If, in univariate histopolation \cite{BE23}, this can be handled by identifying some relevant subclasses of segments, in the multivariate one even this simplification fails. As a consequence, in Section \ref{sect:numerics} we analyze the counterpart of nodal techniques for the identification of approximated Leja sequences \cite{BDSV10} and Fekete sets \cite{BCLSV11}. These algorithmic techniques, partially motivated by stability results \cite{PiazzonVianello}, yield in fact a collection of supports with a low Lebesgue constant.

\begin{remark}
To make the description more readable, the dependence of the functions and the polynomials on the variable $ x $ (which also denotes multivariate coordinates $ x_1, \ldots, x_n $) will be frequently suppressed.
\end{remark}

\section{Polynomial approximation} \label{sect:uniform}

Uniform approximation \cite[Chapter 6]{Davis} aims at controlling the sup-norm of the interpolation procedure, that is, if $ p $ is the polynomial that satisfies \eqref{eq:areamatching}, it deals with the quantity $ \Vert f - p \Vert_{\infty} $. This is generally performed via the analysis of two inequalities. The first one controls the stability of the problem: if the data collected are not exact, the perturbation $ \Vert f - \widetilde{f} \Vert_{\infty} < \varepsilon $ propagates proportionally to the norm of the interpolation operator
\begin{equation} \label{eq:opnorm}
\Vert \Pi \Vert_{\mathrm{op}} := \sup_{f} \frac{\Vert \Pi f \Vert_\infty}{\Vert f \Vert_\infty} = \sup_{\Vert f \Vert_\infty = 1} \Vert \Pi f \Vert_\infty .
\end{equation}
Indeed, one has
\begin{equation} \label{eq:stab}
    \Vert \Pi f - \Pi \widetilde{f} \Vert_{\infty} = \Vert \Pi (f - \widetilde{f}) \Vert_{\infty} \leq \Vert \Pi \Vert_{\mathrm{op}} \Vert f - \widetilde{f} \Vert_{\infty} .
\end{equation}
Further, if $ \Pi $ is a projector, one has that
\begin{equation}  \label{eq:approx}
    \Vert f - \Pi f \Vert_{\infty} \leq \Vert f - p \Vert_{\infty} + \Vert \Pi f - \Pi p \Vert_{\infty} \leq (1+ \Vert \Pi \Vert_{\mathrm{op}} ) \inf_{p \in \mathbb{P}_d (\Omega)} \Vert f - p \Vert_\infty .
\end{equation}
This is the case of the interpolator defined by the set of conditions in Eq. \eqref{eq:areamatching}.

\begin{proposition} \label{prop:projector}
    Suppose $ \{ K_i \}_{i=1}^N $ is unisolvent for $ \P_d (\Omega) $. The interpolator $ \Pi: C^{0} (\Omega) \to \P_d (\Omega) $ defined in \eqref{eq:Lagrangerepresentation} is a projector.
\end{proposition}

\begin{proof}
    By plugging \eqref{eq:Lagrangeduality} into \eqref{eq:Lagrangerepresentation}, one immediately retrieves that $ \Pi \ell_{K_i} (x) = \ell_{K_i} (x) $. We then compute
    \begin{align*}
        \Pi \left(\Pi f\right) & = \sum_{j=1}^N \int_{K_j} \left( \sum_{i=1}^N \left( \int_{K_i} f \right) \ell_{K_i} \right) \ell_{K_j} = \sum_{j=1}^N \left( \sum_{i=1}^N \left( \int_{K_i} f \right)  \int_{K_j} \ell_{K_i} \right) \ell_{K_j} \\ & = \sum_{j=1}^N \left( \int_{K_j} f \right) \ell_{K_j} = \Pi f,
    \end{align*}
    which shows the claim. \qed
\end{proof}

The inequalities \eqref{eq:stab} and \eqref{eq:approx}, together with the representation \eqref{eq:Lagrangerepresentation}, evidence the role of the operator norm $ \Vert \Pi \Vert_{\mathrm{op}} $ as the parameter that controls the numerical conditioning of the interpolation procedure. 

In nodal Lagrange interpolation, its norm $ \Vert \Pi \Vert_{\mathrm{op}} $ may be expressed in terms of the Lagrange cardinal functions $ \{ \ell_{\xi_i} \}_{i=1}^N $ as $ \Vert \Pi \Vert_{\mathrm{op}} = \sup_{x\in\Omega} \sum_{i=1}^N | \ell_{\xi_i} (x) | $, provided that the norm \eqref{eq:equivalencenorm} is considered. This quantity is highly sensitive to the choice of the nodes $ \xi_i $'s and it is used to detect well-conditioned sets of nodes.

\subsection{The Lebesgue constant}

To select good sets of supports for diffused data, we adapt the quantity proposed in \cite{ARR20} as
\begin{equation} \label{eq:Leb}
    \Lambda_d := \sup_{K \subset \Omega} \frac{1}{\mu(K)} \sum_{i=1}^N \mu(K_i) \left| \int_{K} \ell_{K_i} \right| = \sup_{\xi \in \Omega} \sum_{i=1}^N \mu(K_i) \left| \ell_{K_i} (\xi) \right|,
\end{equation}
the last equality being \eqref{eq:equivalencenorm} and the 
Lagrange functions $ \ell_{K_i} $ are defined via the duality relation \eqref{eq:Lagrangeduality}.

Apart from few situations, mostly one-dimensional (see e.g. \cite[Theorem 5.2]{Robidoux} or \cite[Proposition 3.6]{BE23}), explicit expressions for histopolation-based Lagrange functions $ \ell_{K_i} $ are not available. They may nevertheless be recovered by selecting any convenient basis $ \{ p_j \}_{j=1}^N $ for $ \P_d (\R^n) $ and constructing the \emph{Vandermonde matrix}
\begin{equation} \label{eq:VdM}
    \boldsymbol V_{i,j} := \int_{K_i} p_j ,
\end{equation}
so that $ \ell_{K_i} (x) = \sum_{j=1}^N \boldsymbol V^{-T}_{i,j} p_j (x) $. Note that this requires the invertibility of the Vandermonde matrix, which is independent of the basis chosen for $ \P_d (\R^n) $ and is in fact an equivalent condition for unisolvence, see e.g. \cite[Lemma 3.2.2]{Atkinson}.

\begin{remark} \label{rmk:normalizedbasis}
    The term $ \mu(K_i) $ in \eqref{eq:Leb} may be incorporated in the Lagrange functions modifying \eqref{eq:Lagrangeduality} as 
    $$ \frac{1}{\mu(K_i)} \int_{K_i} \widetilde{\ell}_{K_j} = \delta_{i,j} ,$$
    obtaining the equivalent expression $ \Lambda_d = \sup_{\xi \in \Omega} \sum_{i=1}^N |\widetilde{\ell}_{K_i} (\xi) | $. In this case, the corresponding Vandermonde matrix scales that defined in Eq. \eqref{eq:VdM} by the diagonal matrix $ \mathrm{diag} (\mu(K_1), \ldots, \mu(K_N)) $.
\end{remark}

Both the characterizations in \eqref{eq:Leb} are interesting: the first one will be used to relate $ \Lambda_d $ with $ \Vert \Pi \Vert_{\mathrm{op}}$, whereas the second one \lbb{helps in reducing} the computational cost in simulations.

\begin{remark} \label{rmk:Ldm1}
    By construction, $ \Lambda_d \geq 1 $ for each $ d $. Indeed, since each $ K_i $ is compact, considering $ K = K_j $ and exploiting the duality relationship of the Lagrange basis, one has
    $$ \Lambda_d = \sup_{K \subset \Omega} \frac{1}{\mu(K)} \sum_{i=1}^N \mu(K_i) \left| \int_{K} \ell_{K_i} \right| \geq \frac{1}{\mu(K_j)} \sum_{i=1}^N \mu(K_i) \left| \int_{K_j} \ell_{K_i} \right| = \frac{\mu(K_j)}{\mu(K_j)} = 1 . $$
\end{remark}

The next result shows that $ \Lambda_d $ indeed overestimates $ \Vert \Pi \Vert_{\mathrm{op}} $.

\begin{proposition} \label{prop:operatornormbound}
    Let $ \{ K_i \}_{i=1}^N $ (possibly not disjoint) be unisolvent for $ \P_d(\R^n) $. Then
    $$ \Vert \Pi \Vert_{\mathrm{op}} \leq \Lambda_d .$$
\end{proposition}

\begin{proof}
    Expanding the definition of $ \Pi $ and the right hand side equality in \eqref{eq:Leb}, we compute
    \begin{align*}
        \Vert \Pi \Vert_{\mathrm{op}} & = \sup_{\Vert f \Vert = 1} \Vert \Pi f \Vert_{\infty} = \sup_{\Vert f \Vert = 1} \sup_{\xi \in \Omega} \left| \sum_{i=1}^N \left( \int_{K_i} f \right) \ell_{K_i} (\xi) \right| \\ 
        & \leq \sup_{\Vert f \Vert = 1} \sup_{K \subset \Omega} \sum_{i=1}^N \left| \int_{K_i} f \right| \left| \int_K \ell_{K_i} \right| \\ & = \sup_{\Vert f \Vert = 1} \sup_{K \subset \Omega} \sum_{i=1}^N \frac{\mu(K_i)}{\mu(K_i)} \left| \int_{K_i} f \right| \left| \int_K \ell_{K_i} \right| \\ &
        \leq \sup_{K \subset \Omega} \sum_{i=1}^N \mu(K_i) \left| \int_K \ell_{K_i} \right| = \Lambda_d .
    \end{align*}
    The claim is proved. \qed
\end{proof}

The numerical experiment depicted in \cite[Fig. 3]{BE23} shows that the distance between $ \Vert \Pi \Vert_{\mathrm{op}} $ and $ \Lambda_d $ increases as the overlap of the supports enlarges. Imposing disjointedness of the supports, we also get the converse inequality.

\begin{theorem} \label{thm:normLeb}
    Let the set $ \{ K_i \}_{i=1}^N $ be unisolvent for $ \P_d (\R^n) $. If \lbb{$ K_i \cap K_j = \emptyset $} for $ i \ne j $, then
    \begin{equation} \label{eq:newLeb}
        \Vert \Pi \Vert_{\mathrm{op}} = \Lambda_d.
    \end{equation}
\end{theorem}

\begin{proof}
Thanks to Proposition \ref{prop:operatornormbound}, we are only left to prove that $ \Vert \Pi \Vert_{\mathrm{op}} \geq \Lambda_d $. To do so, for each $ 0 < \varepsilon < \Lambda_d $ we exhibit a function $ f_\varepsilon \in C^0 (\Omega) $ such that $ \Vert \Pi f_\varepsilon \Vert_{\infty} \geq \Lambda_d - \varepsilon $. In view of Remark \ref{rmk:Ldm1}, for each $ K_i $ we may find an open set $ U_i \subset K_i $ such that $ \mu( K_i \setminus U_i ) \leq \frac{\varepsilon}{\Lambda_d} \mu(K_i) $, or equivalently, $ \mu(U_i) = (1-\frac{\varepsilon}{\Lambda_d}) \mu(K_i) $. Notice that such a construction is well-posed due to the assumptions of Remark \ref{rmk:assumptionscompacts}. Using mollifiers, we may define collections of continuous functions
$$
    f^{(i)}_\varepsilon (x) = \begin{cases}
        \pm 1 \quad & x \in U_i \\
        0 \quad & x \not\in K_i
    \end{cases}
$$
such that $ \Vert f^{(i)}_\varepsilon \Vert_{\infty} = 1 $. Define $ f_\varepsilon := \sum_{i=1}^N f_\varepsilon^{(i)}$. Independently of the choice of the sign made on each $ f_\varepsilon^{(i)} $, the function $ f_\varepsilon $ is continuous since the $K_i$'s are disjoint and the $ f_\varepsilon^{(i)} $'s are continuous. It it also clear that, by construction, $ \Vert f_\varepsilon \Vert_{\infty} = 1 $. Further, for any $ K \subset \Omega $, there exists a function $ f_\varepsilon $ that satisfies
    $$ \mathrm{sgn} \int_{K_i}f_\varepsilon = \mathrm{sgn} \int_{K } \ell_{K_i} $$
    for each $ i = 1, \ldots, N $ (the case $ \int_K \ell_{K_i} = 0 $ can be arbitrarily treated). By construction, one also has
    $$ \left\vert \int_{K_i} f_\varepsilon \right\vert = \left\vert \int_{K_i} f_\varepsilon^{(i)} \right\vert \geq \mu(U_i) = \left( 1-\frac{\varepsilon}{\Lambda_d} \right) \mu(K_i) .$$
    Using the above two facts and exploiting the equality in \eqref{eq:Leb}, we compute
    \begin{align*}
        \Vert \Pi \Vert_{\mathrm{op}} \geq \Vert \Pi f_\varepsilon \Vert_{\infty} & = \sup_{K \in \subset \Omega} \frac{1}{\mu(K)} \left\vert \int_K \sum_{i=1}^N \left( \int_{K_i} f_\varepsilon \right) \ell_{K_i} \right\vert \\
        & = \sup_{K \subset \Omega } \frac{1}{\mu(K)} \left\vert \sum_{i=1}^N \left( \int_{K_i} f_\varepsilon \right) \int_K \ell_{K_i} \right\vert  \\  
        &\geq \sup_{K \subset \Omega } \frac{1}{\mu(K)} \left\vert \sum_{i=1}^N \mu(U_i) | \int_K \ell_{K_i} \right\vert \\
        & = \sup_{K \subset \Omega } \frac{1}{\mu(K)} \left\vert \sum_{i=1}^N \left( 1-\frac{\varepsilon}{\Lambda_d} \right) \mu(K_i) \int_K \ell_{K_i} \right\vert \\
        & = \left( 1-\frac{\varepsilon}{\Lambda_d} \right) \sup_{K \subset \Omega } \frac{1}{\mu(K)} \left\vert \sum_{i=1}^N  \mu(K_i) \int_K \ell_{K_i} \right\vert = \Lambda_d - \varepsilon .
    \end{align*}
The claim is proved.  \qed
\end{proof}

\subsection{Invariance under affine transformations}

Following the discussion in \cite[Section 6]{ABR20}, one sees that Lebesgue constants associated with non point-wise objects may not be invariant under the choice of the reference domain. In thi section, we show that the quantity \eqref{eq:Leb} depends only on the positioning of $ K_1, \ldots, K_N $ inside $ \Omega $ and not on the placement of $ \Omega $ in $ \R^n $. 
This follows from the fact that if $ \Omega $ is rigidly mapped to another reference domain by $ \varphi $, the Lebesgue constant associated with $ K_1, \ldots, K_N $ in $ \Omega $ coincides with the Lebesgue constant associated with $ \varphi(K_1), \ldots, \varphi(K_N) $ in $\widehat{\Omega} := \varphi (\Omega) $.

\begin{proposition} \label{prop:invarianceLeb}
    Let $ \varphi: \Omega \to \widehat{\Omega} $ be an affine transformation $ \varphi(x) = Ax + b $ such that $ \det A \ne 0 $. Let
    $$ \Lambda_d^{\varphi} := \sup_{\widehat{K} \in \widehat{\Omega}}  \frac{1}{\mu(\widehat{K})} \sum_{i=1}^N \mu(\varphi(K_i)) \left \vert \int_{\widehat{K}} \ell_{\varphi(K_i)} \right \vert $$
    be the quantity \eqref{eq:newLeb} associated with $ \left\{ \varphi (K_i) \right\}_{i=1}^N $. Then
    $$ \Lambda_d = \Lambda_d^{\varphi} .$$
\end{proposition}

\begin{proof}
By the change of variable, we have
$$ \delta_{i,j} = \int_{K_j} \ell_{K_i} = \int_{\varphi(K_j)} \frac{1}{|\det A|} \ell_{K_i} \circ \varphi^{-1} .$$
Put $ \ell_{\varphi(K_i)} := | \det A |^{-1} \ell_{K_i} \circ \varphi^{-1} $. Since $ \varphi^{-1} $ is a non-degenerate affine map, $ \left\{ \ell_{\varphi(K_i)} \right\}_{i=1}^N $ is a basis for polynomials of degree $ d $ satisfying the duality constraint \eqref{eq:Lagrangeduality}. As a consequence, it is the Lagrange basis for $ \P_d (\widehat{\Omega}) $ associated with $ \left\{\varphi(K_i) \right\}_{i=1}^N$. By the same calculation, we also get that, for each $ K \subset \Omega $,
    $$ \int_{K} \ell_{K_i} = \int_{\varphi(K)} \ell_{\varphi(K_i)}, \qquad i = 1, \ldots, N .$$
    To conclude the proof, it is now sufficient to write $ \widehat{K} = \varphi(K) $ and observe that $ \mu(\varphi(K)) = | \det A | \mu (K) $ for each $ K \subset \Omega $. Applying the above equalities and the definition in Eq. \eqref{eq:Leb}, we readily compute
    \begin{align*}
        \Lambda_d^\varphi & := \sup_{\widehat{K} \subset \widehat{\Omega}} \frac{1}{\mu (\widehat{K})} \sum_{i=1}^N \mu(\varphi(K_i)) \left| \int_{\widehat{K}} \ell_{\varphi(K_i)} \right| \\
        & = \sup_{\varphi(K) \subset \widehat{\Omega}} \frac{1}{\mu(\varphi(K))} \sum_{i=1}^N \mu(\varphi(K_i)) \left| \int_{\varphi(K)} \ell_{\varphi(K_i)} \right| \\
        & = \sup_{K \subset \Omega} \frac{1}{|\det A| \mu(K)} \sum_{i=1}^N |\det A| \mu(K_i) \left| \int_K \ell_{K_i} \right| \\
        & = \sup_{K \subset \Omega} \frac{1}{\mu(K)} \sum_{i=1}^N \mu(K_i) \left| \int_K \ell_{K_i} \right| = \Lambda_d,
    \end{align*}
    and the proof is concluded.  \qed
\end{proof}


\subsection{Numerical estimation of the Lebesgue constant}

So far, we have dealt with the definition of Lebesgue constants on the continuous level. In this section we provide an implementation strategy. Due to the lack of general closed formulae for the Lagrange basis, this technique requires the inversion of the Vandermonde matrix, which can rapidly become ill-conditioned \cite{SerraConditioning}.

\subsubsection{An implementation strategy} \label{sect:EstimationLeb}

To estimate Lebesgue constants associated with some unisolvent sets $ \{ K_i \}_{i=1}^N $ on a reference domain $ \Omega $ we discretize Eq. \eqref{eq:Leb}. For what concerns the right hand side, this consists in selecting a large but finite family of supports $ \mathcal{S} := \{ S_\ell \}_{\ell=1}^M $ with $ M \gg N $. Let us expand the leftmost equality in \eqref{eq:Leb} as
\begin{align*}
    \Lambda_d & \approx \sup_{S_\ell \in \mathcal{S}} \frac{1}{\mu({S_\ell})} \sum_{i=1}^N \mu(K_i) \left\vert \int_{S_\ell} \ell_{K_i} \right\vert = \sup_{S_\ell \in \mathcal{S}} \frac{1}{\mu({S_\ell})} \sum_{i=1}^N \mu(K_i) \left\vert \int_{S_\ell} \sum_{j=1}^n \boldsymbol V^{-T}_{i,j} p_j \right\vert \\
    & = \sup_{S_\ell \in \mathcal{S}} \frac{1}{\mu({S_\ell})} \sum_{i=1}^N \mu(K_i) \left\vert \sum_{j=1}^n \boldsymbol V^{-T}_{i,j} \int_{S_\ell}  p_j \right\vert,
\end{align*}
where $ \{ p_j \}_{j=1}^N $ is any basis for $ \P_d (\Omega) $. We will discuss in the next Section \ref{sect:choiceofthebasis} a convenient choice. Define the $ M \times N $ matrix
$$ \boldsymbol W_{\ell,j}^{\mathcal{S}} := \int_{S_\ell} p_j $$
and the square diagonal matrices
$$ \boldsymbol{S} = \begin{pmatrix} \mu(S_1)^{-1} & 0 & \ldots & 0 \\
0 & \ddots & \ddots & 0 \\
0 & \ddots & \ddots & 0 \\
0 & \ldots & 0 & \mu(S_M)^{-1} 
\end{pmatrix} \quad \text{and} \quad
\boldsymbol{K} = \begin{pmatrix} \mu(K_1) & 0 & \ldots & 0 \\
0 & \ddots & \ddots & 0 \\
0 & \ddots & \ddots & 0 \\
0 & \ldots & 0 & \mu(K_N)
\end{pmatrix} .
$$
It follows that, for each $ \ell $,
$$ \frac{1}{\mu({S_\ell})} \sum_{i=1}^N \mu(K_i) \left\vert \int_{S_\ell} \ell_{K_i} \right\vert = \sum_{t=1}^N \left\vert \left(\boldsymbol{S} \boldsymbol W^{\mathcal{S}} \boldsymbol V^{-1} {\boldsymbol{K}}\right)_{\ell,t} \right\vert ,$$
whence
\begin{equation} \label{eq:longLambda}
\Lambda_d \approx \Vert \boldsymbol{S} \boldsymbol W^{\mathcal{S}} \boldsymbol V^{-1} {\boldsymbol{K}} \Vert_{\lbb{\infty}}, 
\end{equation}
being $ \Vert \cdot \Vert_{\lbb{\infty}} $ the $\lbb{\infty}$-norm for matrices (i.e. the maximum by row-wise sum of absolute values, see e.g. \cite[Eq. (2.3.10)]{Golub}). Even if one has at disposal and uses an exact quadrature formula (as, e.g., those used in \cite{AESV24,ASV20}) for the elements of $ \mathcal{S} $, so that no approximation error is carried on $ \boldsymbol W^{\mathcal{S}} $, it is clear that the assembly of such a matrix is expensive. We hence discretize also the rightmost term of \eqref{eq:Leb} by considering a comparably large set of nodes $ \mathcal{X} := \{ \xi_i \}_{i=1}^M $, with $ M \gg N $. By the very same reasoning, we obtain that
\begin{align*}
    \Lambda_d & \approx \sup_{\xi_\ell \in \mathcal{X}} \sum_{i=1}^N \mu(K_i) \left\vert \ell_{K_i} (\xi_\ell) \right\vert = \sup_{\xi_\ell \in \mathcal{X}} \sum_{i=1}^N \mu(K_i) \left\vert \sum_{j=1}^n \boldsymbol V^{-T}_{i,j} p_j (\xi_\ell) \right\vert  \\
    & = \sup_{\xi_\ell \in \mathcal{X}} \sum_{i=1}^N \mu(K_i) \left\vert \sum_{j=1}^n \boldsymbol V^{-T}_{i,j}  p_j (\xi_\ell) \right\vert.
\end{align*}
In this case, defining the $ M \times N $ matrix 
$$ \boldsymbol W_{\ell,j}^{\mathcal{X}} := p_j (\xi_\ell) $$
one readily obtains
\begin{equation} \label{eq:shortLambda}
\Lambda_d \approx \Vert \boldsymbol W^{\mathcal{X}} \boldsymbol V^{-1} \boldsymbol{K} \Vert_{\lbb{\infty}} .
\end{equation}
Comparing \eqref{eq:shortLambda} with \eqref{eq:longLambda}, it is evident that the latter formulation is sensibly cheaper, as it avoids the computation of several integrals and does not require the first matrix-matrix product. Since the equalities in \eqref{eq:Leb} hold only if all compacts $K$ and nodes $ \xi $ in $ \Omega $ are considered, \eqref{eq:shortLambda} and \eqref{eq:longLambda} may differ. Numerical verifications show that the distance between such quantities is in fact negligible.

\subsubsection{A remark on the conditioning of the Vandermonde matrix} \label{sect:choiceofthebasis}

Obstructions in the computation of Lagrange bases for histopolation are well-known \cite{BE23,Gerritsma}. 
%
Due to the lack of closed formulae, one is forced to fix a convenient basis $ \{ p_i \}_{i=1}^N $ for $ \P_d (\Omega)$, compute the relative Vandermonde matrix \eqref{eq:VdM} and invert it to obtain the coefficients of \eqref{eq:Lagrangeduality}, as described after Eq. \eqref{eq:VdM}. Of course, the inversion of an ill-conditioned matrix may cause a severe instability. This topic is widely studied in the literature; in particular, it has recently been observed that up to a mild degree this effect is not particularly dangerous \cite{MonomialConditioning}.

To avoid the propagation of errors and instability, among easy-to-write bases, we select the one that shows the smallest conditioning of the resulting Vandermonde matrix. Two natural candidates are thus the monomial $ \mathcal{M}_d := \{ x^\alpha y^\beta \} $ and the Chebyshev $ \mathcal{T}_d := \{ T_\alpha (x) T_\beta (y) \} $ bases, with $ \alpha + \beta \leq d $. 

\begin{figure}[!h]
   \centering
        {\includegraphics[width=5.9cm, clip,trim = 0 0 0 0]{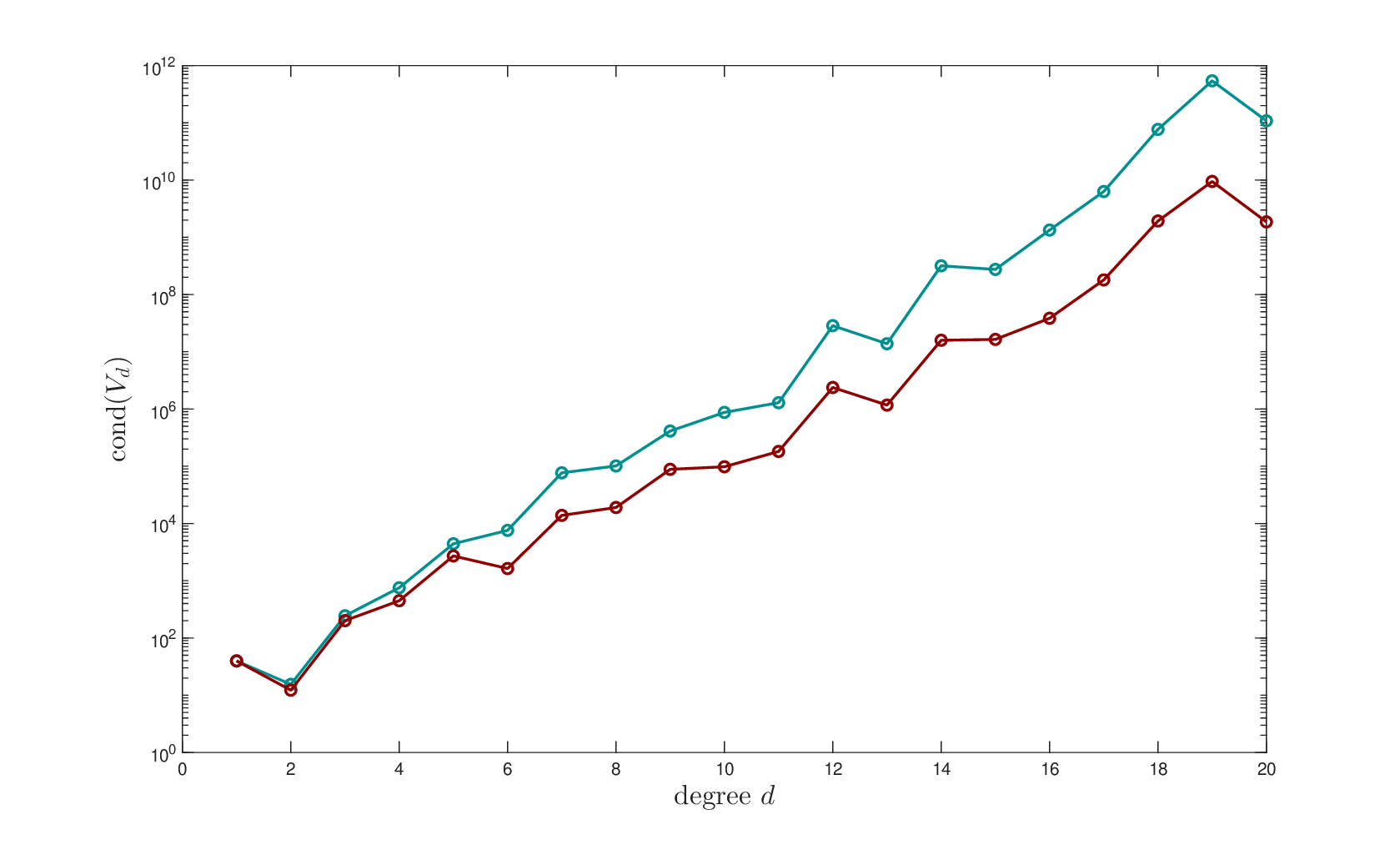}}
     \quad
        {\includegraphics[width=5.9cm, clip,trim = 0 0 0 0]{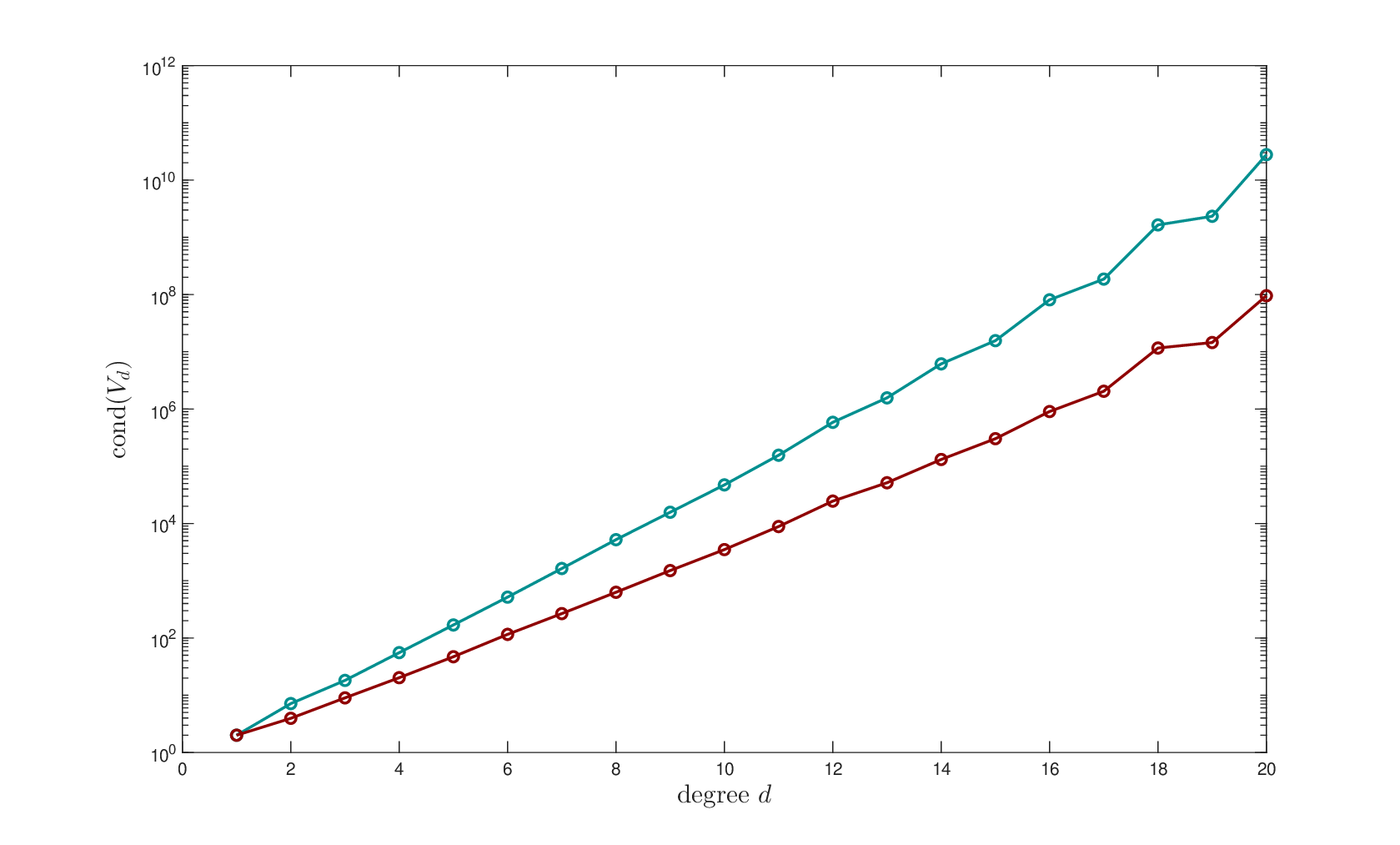}}     
\caption{Conditioning of the Vandermonde matrix: Chebyshev polynomials (red) vs monomials (blue). Supports are discs centered at Halton points (left) and Chebyshev points in the reference disc (right).} \label{fig:cond}
\end{figure}

In Figure \ref{fig:cond} we depict the trend of the conditioning of the Vandermonde matrix $ \boldsymbol V $ with respect to the monomial basis (blue) and the Chebyshev basis, for discs of fixed radius centered at Halton points \cite{Halton} (which are quasi-random points) and points defined in \cite{BX03} (due to their construction, we will refer to such points as Chebyshev orbits). In both cases Chebyshev polynomials lead to a better conditioned Vandermonde matrix; in the numerical section we will thus stick with this choice.

\section{Two approaches to unisolvence} \label{sect:unisolvence}

The detection of unisolvent sets is the starting step of any interpolation problem. In the multidimensional setting, already in the nodal framework, few explicit unisolvent sets of points for the polynomial space $ \P_d (\Omega) $ are known. Nevertheless, it is well known that any collection of nodes is \emph{generically} unisolvent, meaning that the measure of the manifold of all non-unisolvent sets for polynomial interpolation is zero.

\begin{remark} \label{rmk:density}
We already observed that unisolvence is equivalent to the invertibility of the Vandermonde matrix $ \boldsymbol V $ defined in \eqref{eq:VdM}. As a consequence of the mean value theorem, the determinant $ \det : (\R^n)^N \to \R$ is a polynomial map whose variables are the points $ \xi_i \in K_i $ verifying the statement of the mean value theorem. Since any collection of nodes is generically unisolvent, in the above defined sense, it follows that the locus of non-unisolvent collection of compact sets (represented by $ \det^{-1} (0)$) is a Zariski closed subset of $ (\R^n)^N $ and hence has Lebesgue measure zero. 
As a consequence, a generic choice of compacts $ \{ K_i \}_{i=1}^N $ is unisolvent. 
\end{remark}

\subsection{Unisolvence by translation}

Despite the density result of Remark \ref{rmk:density}, the identification of classes of unisolvent supports is a hard task, if not connected with nodal results. Also, such class of results is of limited direct applicability, as the mean value theorem does not allow for a precise localization of the nodes. If we fix the a reference compact set $ K $ and consider its translates, we may relate with nodal interpolation results. To begin with, we extend \cite[Lemma 3.12]{ChR16} to the case of generic compacts. We adopt the usual notation $ K + \xi := \{ x + \xi \mid x \in K \} $.

\begin{lemma} \label{lem:alltranslated}
    Let $ K \subset \R^n $ be a compact set and let $ p \in \P_d (\R^n) $. If
    \begin{equation} \label{eq:genericxi}
        \int_{K + \xi} p(x) \de x = 0
    \end{equation}
    for each $ \xi \in \R^n $, then $ p(x) = 0 $.
\end{lemma}

\begin{proof}
    Suppose, by contradiction, that $ p $ is not the zero polynomial. Let us denote by $ \mathcal{Z}[p] $ the zero locus of $ p $, which is an algebraic variety. Since $p$ is not identical to zero, the measure of $\mathcal{Z}[p]$ is equal to zero. Hence there exists $ \xi \in \mathbb{R}^n $ such that $ (K + \xi) \cap \mathcal{Z}[p] = \emptyset $. As a consequence, $ \mathrm{sgn} \left( p \right) $ is constant in $ K + \xi $ and its interior, which is open and non empty, being the measure of $ K  $ finite and greater than zero. Thus Eq. \eqref{eq:genericxi} implies that $ p (x) = 0 $ for each $ x $ in such an open set and therefore $ p(x) = 0 $  for each $ x \in \mathbb{R}^n $.  \qed
\end{proof}

\begin{proposition} \label{prop:unisolvencepoints}
    Let $ K \subset \R^n $ be a compact set. Adopting the usual convention 
    $$ AK + \xi := \{ Ax + \xi \mid x \in K \}, $$
    suppose that there exists a non-degenerate family of affinities $ \varphi_{i} (x) = A x + \xi_i $ such that $ K_i := \varphi_i (K) = A K + \xi_i $ for $ i = 1, \ldots, N $. If
    \begin{enumerate}
        \item[(i)] $ \{ \xi_i \}_{i=1}^N $ is a unisolvent set for nodal interpolation in $ \P_d (\Omega) $,
    \end{enumerate}
    then
    \begin{enumerate}
        \item[(ii)] $ \{ K_i \}_{i=1}^N $ is a unisolvent set for integral interpolation in $ \P_d (\Omega) $.
    \end{enumerate}
\end{proposition}

\begin{proof}
    We shall prove that 
    $$ \int_{K_i} p(x) \de x = 0 \quad i = 1, \ldots, N \qquad \Longrightarrow \qquad p(x) = 0 \quad \forall x \in \Omega. $$
    Consider the map 
    $$ q: \quad \xi \mapsto \int_{AK + \xi} p(x) \de x . $$
    Since $ \xi \mapsto AK + \xi $ is an affinity, by the change of variable formula one deduces that $ q(\xi) $ is a polynomial map. Further, since $ p \in \P_d $, also $ q \in \P_d $. By construction, $ q $ vanishes at $ \{ \xi_i \}_{i=1}^N $, i.e. $ q(\xi_i) = 0 $ for $ i = 1, \ldots, N $, which is a unisolvent set for $ \P_d $, so $q$ is the null polynomial; namely $ q(\xi) = 0 $ for each $ \xi \in \R^n $. 
    Thus 
    $$ 0 = q(\xi) = \int_{AK + \xi} p(x) \de x \qquad \forall \xi \in \R^n .$$
    Since $ A $ is fixed, the claim follows from Lemma \ref{lem:alltranslated}.  \qed
\end{proof}

\begin{remark}
    In the univariate framework, the converse implication $(ii) \Longrightarrow (i) $ holds. Indeed, if $ \int_{K_i} p(x) \de x = 0 $ implies that $ p(x)= 0 $, in particular all the supports $ K_i $ must be distinct (otherwise we find a repeated condition and the problem becomes underdetermined). Hence $K_i \ne K_j $ if $i \ne j$. Since $K_i =AK + \xi_i $ and $K_j =AK+ \xi_j$, one immediately finds $\xi_i \ne \xi_j $ for $i \ne j $. As a consequence, the set  $ \{ \xi_i \}_{i=1}^N $ contains pairwise distinct points and is therefore unisolvent for nodal interpolation. This partially extends the characterization of unisolvent segments in \cite{BE23}.
\end{remark}

The above result makes it possible to invoke the literature on unisolvent sets for multivariate interpolation, for instance \cite{BSV12,Hesthaven98,MS19}, at the price of requiring that all the supports are the image through $ \varphi $ of the same reference compact set $ K $. The way the compact set $ K $ is selected does not affect unisolvence, provided that it satisfies the assumptions of Remark \ref{rmk:assumptionscompacts}.

\subsection{Unisolvence by algebraic varieties}

The constraint of fixing the reference compact set $ K $ once and for all can be replaced if nodes may be collected into subsets lying on \emph{(affine) algebraic varieties}. To do so, recall that an affine algebraic variety in $ \R^n $ is the zero locus of a finite collection of polynomials in $ n $ variables \cite[p. 3]{Hartshorne}. In what follows, we will consider algebraic varieties $ \mathcal{V} $ defined by just one polynomial $ P \in \P_{\dim \mathcal{V}} (\R^n) $, for short we shall write $ \mathcal{V} := \{ P = 0 \} $. Proposition \ref{prop:unisolvencepoints} immediately yields the following result.

\begin{corollary} \label{cor:factorbit}
    Let $ p \in \P_d (\R^n) $ and let $ \mathcal{V} := \{ P = 0 \} $. If $ \mathcal{X} = \{ \xi_i \}_{i=1}^{\dim \P_d (\mathcal{V})} $ is a nodal unisolvent set for $ \P_d (\mathcal{V}) $ and
    $$ \int_{K_i} p(x) \de x = 0 \qquad i = 1, \ldots, \dim \mathcal{V} ,$$
    $ K_i $ being defined as in Proposition \ref{prop:unisolvencepoints}. Then $ P $ divides $ p $.
\end{corollary}

\begin{proof}
    Since $ \mathcal{X} $ is unisolvent for $ \P_d (\mathcal{V}) $, by Proposition \ref{prop:unisolvencepoints} we have that
    $$ \int_{K_i} p(x) \de x = 0 \qquad i = 1, \ldots, \dim \mathcal{V} ,$$
    implies that $ p (x) = 0 $ for each $ x \in \mathcal{V} $. Since $ \mathcal{V} = \{ P = 0 \} $, there exists $ q \in \P_{d-\dim \mathcal{V}} (\R^n) $ such that $ p = P q $, which shows the claim.  \qed
\end{proof}

Corollary \ref{cor:factorbit} allows to split the global interpolation problem into subproblems onto algebraic varieties, provided that they do not have common components. This means that, given $ \mathcal{V} = \{ P = 0 \} $ and $ \mathcal{W} = \{ Q = 0 \} $, 
there does not exist any polynomial $ R $ such that either $ P = R Q $ or $ Q = R P $.

For each variety $ \mathcal{V}_j $, we select a reference compact $ \widehat{K}_j $ and we translate it by vectors $ \{ \xi_j^i \}_{i=1}^{\dim \P (\mathcal{V}_j) } $. This moves the global unisolvence conditions to problems on the single varieties. We shall thus require that collections $ \mathcal{X}_j := \{ \xi_j^i \}_{i=1}^{\dim \P (\mathcal{V}_j) } $ are unisolvent for the spaces $ \P (\mathcal{V}_j) $ (for some polynomial degree).

Clearly, if $ \mathcal{V}_j = \R^n $, the next result is equivalent to Proposition \ref{prop:unisolvencepoints}. In Remark \ref{rmk:exampler2} we instead give an example in which the splitted interpolation problem is sensibly simpler than the global one.

\begin{theorem} \label{thm:algvar}
    Let $ p \in \P_d (\R^n) $. Let $ \mathcal{V}_1 = \{ P_1 = 0 \}, \ldots, \mathcal{V}_s = \{ P_s = 0 \} $ be affine algebraic varieties with pairwise no common components. Suppose that $ P_j \in \P_{d_j} (\R^n) $ and $ \sum_{j=1}^s d_j > d $. Suppose that, for each $ \mathcal{V}_j $, there exists a collection of nodes $ \mathcal{X}_j := \{ \xi_j^1, \ldots, \xi_{j}^{N_j} \} $ which is unisolvent for $ \P_{d - \sum_{i<j} d_i} (\mathcal{V}_j) $. If, for $ j = 1, \ldots, s $,
    $$ \int_{K_{j}^i} p(x) \de x \ \lbb{ =0} \qquad i = 1, \ldots, N_j := \dim \P_{d - \sum_{i<j} d_i} (\mathcal{V}_j),$$
    being $ K_{j}^i := \varphi_{j}^i (K) = A_j x + \xi_j^i $ with $ \det A_j \ne 0 $, then $ p (x) = 0 $.
\end{theorem}

\begin{proof}
    Up to relabelling, we may reorder the varieties $ \mathcal{V}_i $ so that $ j < k $ if $ d_j > d_k $. Consider $ j = 1 $. Since
    $$ \int_{K_{1}^i} p(x) \de x = 0 \qquad i = 1, \ldots N_1 ,$$
    by Corollary \ref{cor:factorbit} we have that $ P_1 $ divides $ p $; equivalently, we may write $ p = P_1 q $ for some $ q \in P_{d - d_1} $. We may now iterate this process up to $ j = s-1 $, to get that
    $$ p = \left( \prod_{j=1}^{s-1} P_j \right) q_{s-1} .$$
    This formula holds because the varieties $ \mathcal{V}_j $ have no common components, hence at the $j$-th step the polynomial $ P_j $, which divides $ p $ by Corollary \ref{cor:factorbit}, must divide $ q_{j-1} $ and not $ \prod_{k=1}^{j-1} P_k $. For the above reasoning, when $ j = s $,
    $$ \int_{K_s^i} p(x) \de x = 0 \qquad i = 1, \ldots, N_s $$
    implies that $ P_s $ divides $ q_s $. But
    $$ \deg q_s = d - \sum_{j=1}^{s-1} d_j < d_s = \deg P_s, $$
    which implies that $ q_s $ is the zero polynomial, and so is $ p $. This proves the claim.  \qed
    \end{proof}


\begin{remark}
    Theorem \ref{thm:algvar} uses Corollary \ref{cor:factorbit} (hence Proposition \ref{prop:unisolvencepoints}) on each algebraic variety $ \mathcal{V}_j $ independently. This means that, for each variety, one may choose a different reference set $ K $.
\end{remark}

\begin{remark}
    Compact sets $ K_j^i $ (associated with the $j$-th variety) need not intersect $ \mathcal{V}_j $, but may intersect other supports and other varieties.
\end{remark}

\begin{remark} \label{rmk:degenerate}
    The condition $ \sum_{j=1}^s d_j > d $ can be replaced by the equality $ \sum_{j=1}^s d_j = d $ plus one spare information (e.g. an evaluation or another integral) on a support that preserves unisolvence. In such a case, the proof reads exactly the same except for the last step. We give an example of this in Proposition \ref{thm:unisolvenceorbits}.
\end{remark}

\subsection{An application to balls} \label{sect:applicationtoballs}

We make the result of Theorem \ref{thm:algvar} more concrete, replacing algebraic varieties with orbits, and considering balls centered at the origin as reference compact sets. In this case, centers are localized on varieties $ \mathcal{V}_j $ (relating to the nodal results in \cite{Bos91}); in turn, we obtain an easier algorithmic description.

The algebraic variety
\begin{equation} \label{eq:Snorbits}
\mathcal{V} (r) := \left\{ x = (x_1, \ldots, x_n) \in \R^n \ \mid \ \sum_{i=1}^n x_i^2 = r^2 \right\} 
\end{equation}
describes the orbit of a point at distance $ r $ from $ 0 $ under the action of the special orthogonal group $ SO(n) $. Note that $ \mathcal{V}(r) = S_r^{n-1} := r S^{n-1} $. Being orbits, $ \mathcal{V}_j := \mathcal{V} (r_j) $ does not intersect $ \mathcal{V}_k := \mathcal{V} (r_k) $ if $ r_j \ne r_k$. Since the polynomial
\begin{equation} \label{eq:orbit}
    P_j (x_1, \ldots, x_n ) := - r_j^2 + \sum_{i=1}^n x_i^2  
\end{equation} 
describing any $ \mathcal{V}_j $ has degree $2$, the dimension of the space of polynomials on $ \mathcal{V}_j $ is
\begin{equation} \label{eq:dimrestr}
    \dim \mathbb{P}_d \left( \mathcal{V}(r) \right) = \dim  \mathbb{P}_d \left( \R^n \right) - \dim \mathbb{P}_{d-2} \left( \R^n \right) = \binom{d+n}{n} - \binom{d+n-2}{n} ,
\end{equation}
see  \cite[Lemma $2.1$]{Bos91}. As reference compact $ K $, we consider the unit $n$-ball $ B := B(0,1) $ and we take the linear part $ A $ of the affinity $ \varphi_j $ to be a scalar multiple of the identity. Hence, $ B_j^i = \varphi_j^i(B) = 
B(\xi_j^i, r_j) $. Further, since $ B $ is centered at $ 0 $, every $ B_j^i $ is centered at a point in $ \mathcal{V}_j $. Unisolvence for the integral interpolation problem on each variety $ \mathcal{V}_j $ is thus granted if centers $ \xi_j^i $, for $ i = 1, \ldots, N_j $ and each $N_j$ being determined by \eqref{eq:dimrestr}, are placed at unisolvent nodes on the sphere $ S^{n-1}_{r_j} $; 
for an account of such nodes see, for instance, \cite{Phung21,Xu03}. Corollary \ref{cor:factorbit} may thus be simplified as follows.

\begin{lemma} \label{lem:factorisationSr} 
Let $ p \in \P_d (\R^n) $ and let $ \mathcal{X} = \{ \xi_i \}_{i=1}^N $ be unisolvent for nodal interpolation in $ \mathbb{P}_d \left( \mathcal{V}(r) \right) $. If, for some $ \bar{r} > 0 $,
$$ \int_{B(\xi_i,\bar{r})} p(x) \de x = 0 \qquad \text{ for each } \xi_i \in \mathcal{X}, $$
then there exists $ q \in P_{d-2} (\R^n) $ such that $ p = Pq $.
\end{lemma}

Our next construction follows Algorithm \ref{alg:orbits}, whose well posedness will be a consequence of Proposition \ref{thm:unisolvenceorbits}. The idea consists in selecting an orbit, saturating it with a unisolvent set, constructing balls with radius $ r_1 $, passing to the next orbit, saturating it with a (smaller, in view of Lemma \ref{lem:factorisationSr} and Eq. \eqref{eq:dimrestr}) unisolvent set, constructing balls with radius $ r_2 $, and iterating the procedure.

\begin{algorithm}
\caption{Construction of the unisolvent set by orbits}\label{alg:orbits}
\begin{algorithmic}
\Require The total degree $ d $
\Ensure A collection of unisolvent balls
\State $ {\bf r} = [ ] $ \Comment{Vector of the radii}
\While{$d \geq 1 $}
\State choose $ r \not\in {\bf r} $ and add $ r $ to $ {\bf r }$ \Comment{Guarantees disjointess of orbits}
\State fix $ r_j $ \Comment{Radius of the balls on this orbit}
\State compute $ d_j := \dim \P_d (\mathcal{V}_j) $ by \eqref{eq:dimrestr} \Comment{Number of balls to place on the $j$-th orbit}
\State choose $c_j$ which are $ d_j $ unisolvent points on the orbit $ \mathcal{V}_j $
\State construct $ d_j $ balls of radius $ r_j $ centered on $ c_j $
\State $ d = d-2 $ \Comment{Every time reduces by $ 2 $ the total degree}
\EndWhile
\If{$ d = 0$} \Comment{Only the degenerate orbit is left}
\State construct a ball of any radius centered at 0
\EndIf 
\end{algorithmic}
\end{algorithm}
Note that, if $ d $ is the total degree of the polynomial $ p $, unisolvence on each variety implies the convergence of the loop in $ d^* := \lceil \frac{d}{2} + 1 \rceil $ iterations, regardless of the dimension of the ambient space. An iterative application of Lemma \ref{lem:factorisationSr} for $ d_j = 0, 2, \ldots, d $ if $ d $ is even or $ d_j = 1, 3, \ldots, d $ if $ d $ is odd is in fact the core of the proof of the next result.
%

\begin{proposition} \label{thm:unisolvenceorbits} Let $ p \in \P_d (\R^n) $. Let $ \mathcal{X} = \{ \xi^i \}_{i=1}^{\dim \P_d (\R^n)} $ be a collection of points such that $ \mathcal{X}_j := \mathcal{V}_j \cap \mathcal{X} $ is unisolvent for $ \P_{d_j} ( \mathcal{V}_j ) $. 
For each $ j $, fix $ r_j > 0 $ and put $ B_j^i := B (\xi_j^i, r_j) $. If
    $$ \int_{B_{j}^i} p(x) \de x = 0 \qquad i = 1, \ldots, N_j := \dim \P_{d_j} ( \mathcal{V}_j ), $$
    for each $ j = 1, \ldots, d^* $, then $ p(x) = 0 $.
\end{proposition}

\begin{proof}
    For $ j = 1, \ldots, d^* $, order orbits $ \mathcal{V}_j $ so that $ j < k $ if $ \# \left( \mathcal{X}_j \right) > \# \left( \mathcal{X}_k \right) $. Consider $ j = 1 $. Since 
    $$ \int_{B_1^i} p(x) \de x = 0 \qquad \text{for each } \xi_1^i \in \mathcal{X}_1 ,$$
    and $ \mathcal{X}_1 $  is a unisolvent set for nodal interpolation in $ \P_d (\mathcal{V}_1) $, by Lemma \ref{lem:factorisationSr} there exists a degree $ 2 $ polynomial $P_1$ as in \eqref{eq:orbit} such that $ p = P_1 q_1 $, whence $ \deg q_1 = d-2 $. 
    Iterating this up to $ j = d^* -1 $, we get
    \begin{equation} \label{eq:fact}
    p  = \left( \prod_{j=1}^{d^* - 1} P_j  \right) q_{d^* -1} .
    \end{equation}
    Since varieties $ \mathcal{V}_j $ are the orbits of the same group action, they do not intersect and thus have no common components. Hence, for $ j = 1, \ldots, d^* -1 $, $ P_j $ does not divide $ \prod_{k=1}^\ell P_k $ for any $ \ell < j $ and thus $ q_{d^* -1} $ is a polynomial of degree $ \bar{d} \leq 1 $. With respect to to the last orbit $ \mathcal{V}_{d^*} $, we have
    $$ p \bigl|_{\mathcal{V}_{d^*}} = c q_{d^*-1}  ,$$
    with $ c = \prod_{j=1}^{d^*-1} (r_j - r_{d^*})^2 $ being a constant. We show, separating the cases $ \bar{d} = 0, 1 $, that     
    $$ \int_{B_{d^*}^i} q_{d^*-1}(x) \de x = 0 \qquad \text{for all } \xi_{d^*}^i \in \mathcal{X}_{d^*} $$
    implies that $ q_{d^*-1} = 0 $. This will complete the proof.
    
    Let us first consider the case $ \bar{d} = 0 $. The polynomial $ q_{d^*-1} $ is therefore a constant, say $ c' $, 
    with vanishing integral on a full measure set. Hence $ c' = 0 $ and, by Eq. \eqref{eq:fact}, $ p (x) = 0 $. 
    
    If we assume $ \bar{d} = 1 $ the claim is a consequence of Lemma \ref{lem:factorisationSr}: the algebraic variety $ \mathcal{V}_{d^*} $ described by $ P_{d^*} $ has degree $ 2 $ and, since it divides $ q_{d^*-1} $, it must be $ q_{d^*-1} = 0 $. Again by Eq. \eqref{eq:fact}, $ p(x) = 0 $.  \qed
\end{proof}

\begin{remark}
    Describing varieties as orbits 
    clarifies that the degenerate orbit in \eqref{eq:orbit} for $ r = 0 $ preserves unisolvence. We stress that, as a consequence of Proposition \ref{thm:unisolvenceorbits}, this point can in fact be any point not lying on any of the $ \mathcal{V}_j $.
\end{remark}

\begin{remark} \label{rmk:exampler2}
    The condition on the unisolvence on the algebraic varieties may significantly simplify the problem. For instance, if $ n = 2 $, it is sufficient to place distinct nodes on each of the circumferences $ \mathcal{V}_j $, with $ j = 1, \ldots, d^* $. The number of nodes to be placed on each $ \mathcal{V}_j $ is further prescribed by \eqref{eq:dimrestr}.
\end{remark}

\section{Numerical experiments} \label{sect:numerics}

To carry out our numerical experiments, we consider $ \Omega \subset \R^2 $ as the unit disc. Although this is not essential under the theoretical point of view of Section \ref{sect:unisolvence}, such a setting matches a well-known topic in approximation theory \cite{MS19,XuChapter}, and also allows to recover some estimates, such as those in \cite{Nevskii,Nevskii2}. Further, since we showed in Proposition \ref{prop:projector} that $ \Pi $ is a projector onto polynomials of degree $ d $, the norm $ \Vert \Pi \Vert_{\mathrm{op}}$ is bounded from below by
\begin{equation} \label{eq:lowerbound}
    \Vert \Pi \Vert_{\mathrm{op}} \geq c_n \cdot \begin{cases}
        \log d, &\text{if } n = 1, \\
        d^\frac{n-1}{2}, &\text{if } n > 1,
    \end{cases}
\end{equation}
as established in \cite{Sundermann}.

\subsection{Lebesgue constants}
The numerical scheme for the estimation of the Lebesgue constants has been discussed in Section \ref{sect:EstimationLeb}. To begin our analysis, we evidence that the quantity $ \Lambda_d $ is little affected by the variation of the size of the supports. Something more definite can be claimed when the $K_i$'s are convex sets. 
Following the nodal stability inequalities proved in \cite{PiazzonVianello}, one extends \cite[Prop. 5]{PreprintBDEN} in the following way.

\begin{proposition} \label{proposition:stability}
    Let $ \mathcal{X} := \{ \xi_1, \ldots, \xi_N \} $ be a set of nodes which is unisolvent for interpolation in $ \P_d (\Omega) $, and let $ \Lambda_d^{\mathcal X} $ denote its nodal Lebesgue constant. Let $ \{ K_1, \ldots, K_N \} $ be a collection of compact convex sets such that $ \xi_i \in K_i $ and $ \mu (K_i \cap K_j) = 0 $ if $ i \ne j $. Let $ r_{\max} := \frac{1}{2} \max_i \mathrm{diam} (K_i) $, and denote by $ \Lambda_d^K $ the Lebesgue constant \eqref{eq:Leb} of $ \{ K_1, \ldots, K_N \} $. For each $ \alpha \in (0,1) $ so that $ r_{\max} < 2 \,c\, d^r \Lambda_d^{\mathcal{X}} \,\alpha $, $ c $ being the Markov constant and $ r $ an appropriate constant depending on $ \Omega $, hence 
    \begin{equation} \Lambda_d^K \leq \frac{1}{1-\alpha} \Lambda_d^{\mathcal X} .\end{equation}
\end{proposition}

The proof of Proposition \ref{proposition:stability} is omitted: it can be directly deduced from \cite{PreprintBDEN}, observing that the diameter of a triangle coincides with its longest edge.

\begin{remark} Roughly speaking, the above result states that Lebesgue constant associated with histopolation on the set $ \{ K_1, \ldots, K_N \} $ is close to the nodal Lebesgue constant of a given set of nodes $ \{ \xi_1, \ldots, \xi_N\} $, provided that $ \xi_i \in K_i $ for $ i = 1, \ldots, N $ and the diameter of the $ K_i$'s is sufficiently small. This is a consequence of the mean value theorem and the stability estimate given in \cite[Proposition 1]{PiazzonVianello}, where a discussion on the role of the Markov constant is also addressed. In particular, since in numerical experiments we are assuming that $ \Omega $ is a convex body, one has $ r = 2 $.
\end{remark}

From now on, to ease the description of the supports, we consider $ K_i = B(\xi_i, r_i) $, so that parameters $ \xi_i $ and $ r_i $ assume a neat meaning. 
This also matches some literature, such as \cite{Volchkov}.

\begin{figure}[!h]
   \centering
   {\includegraphics[width=5.9cm]{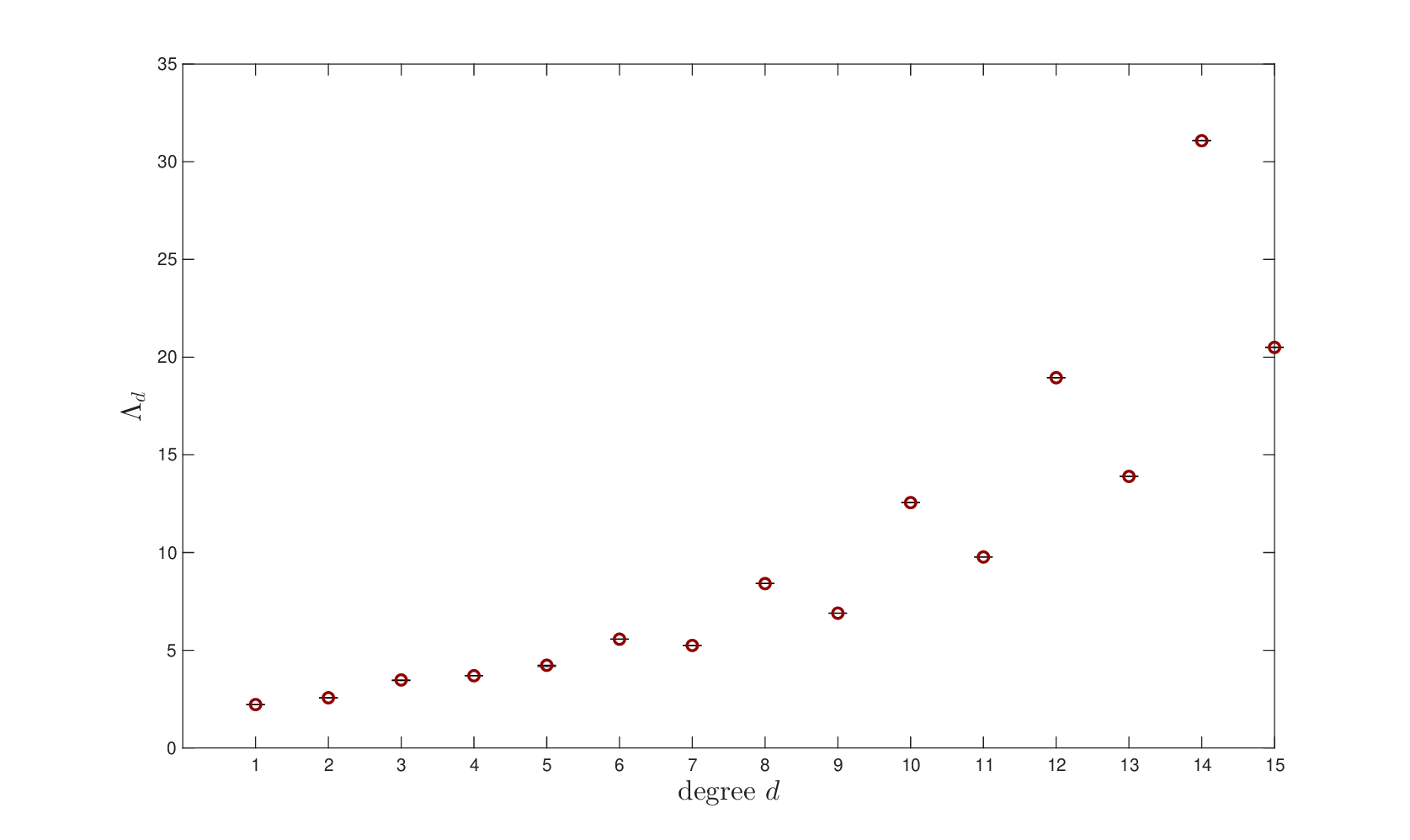}}
   \;
        {\includegraphics[width=5.9cm]{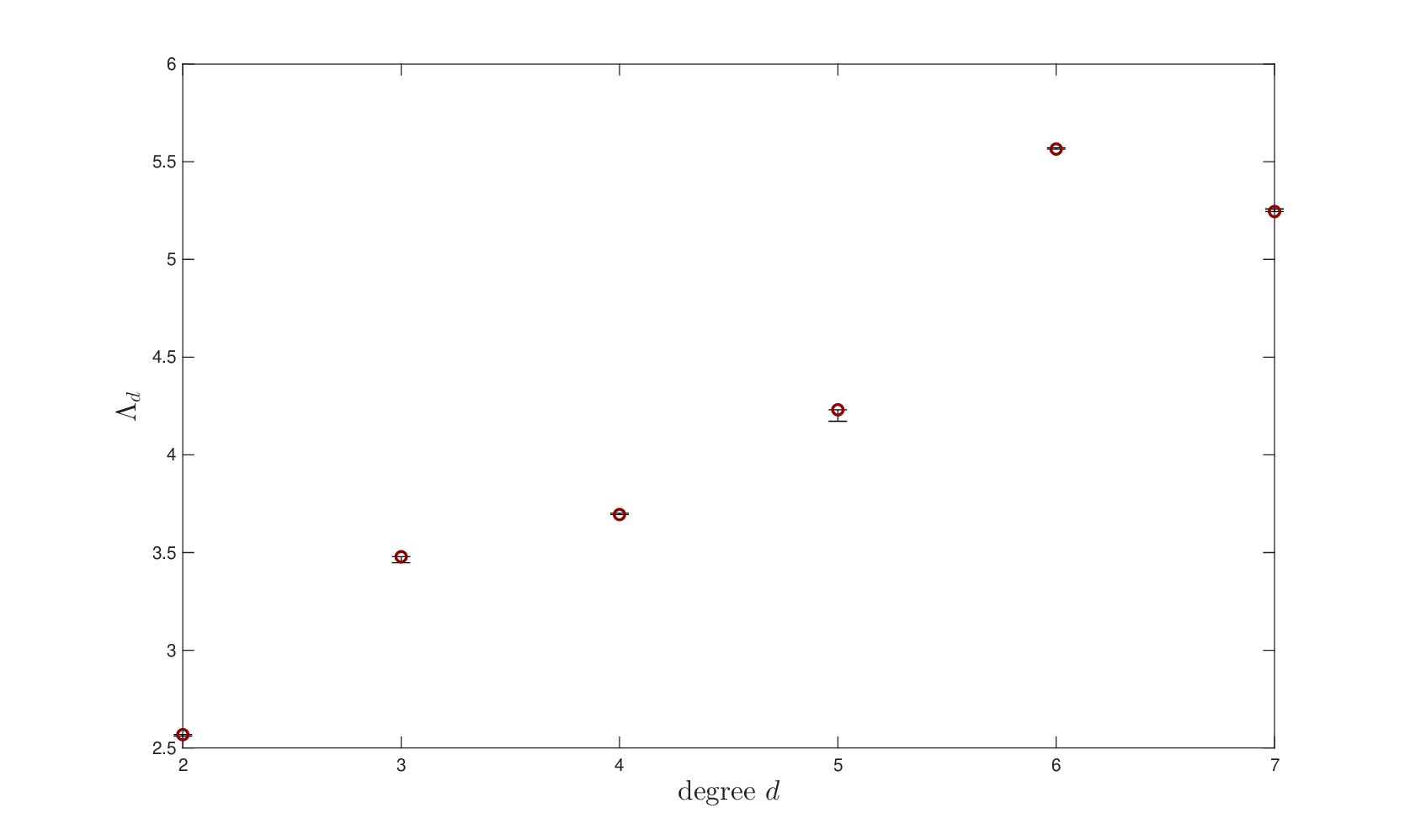}}
\caption{Lebesgue constants of discs centered at Chebyshev points. The nodal Lebesgue constant is represented by the red circles, while the error bars (which are indeed very small) denote the Lebesgue constant of the discs centered at such points, for different radii. The right panel is an enlargement of the left panel. Note that the error bars are almost contained in the circles.} \label{fig:errorbar}
\end{figure}

As a numerical account of Proposition \ref{proposition:stability}, where we consider discs centered at the (unisolvent) nodes of the \textquotedblleft Chebyshev grid\textquotedblright\ proposed by Bojanov and Xu in \cite{BX03}, unisolvence for the histopolation problem associated with such supports can be deduced from both Proposition \ref{prop:unisolvencepoints} and Theorem \ref{thm:algvar}. In particular, any collection of discs of fixed radius $ r > 0 $ centered at $ \mathcal{X} $ is unisolvent for the corresponding histopolation problem. Figure \ref{fig:errorbar} shows that the Lebesgue constant \eqref{eq:Leb}, represented by the error bars (for different radii), is very close to the nodal Lebesgue constant of the set $ \mathcal X$.

\begin{figure}[!h]
   \centering
   {\includegraphics[width=6.2cm]{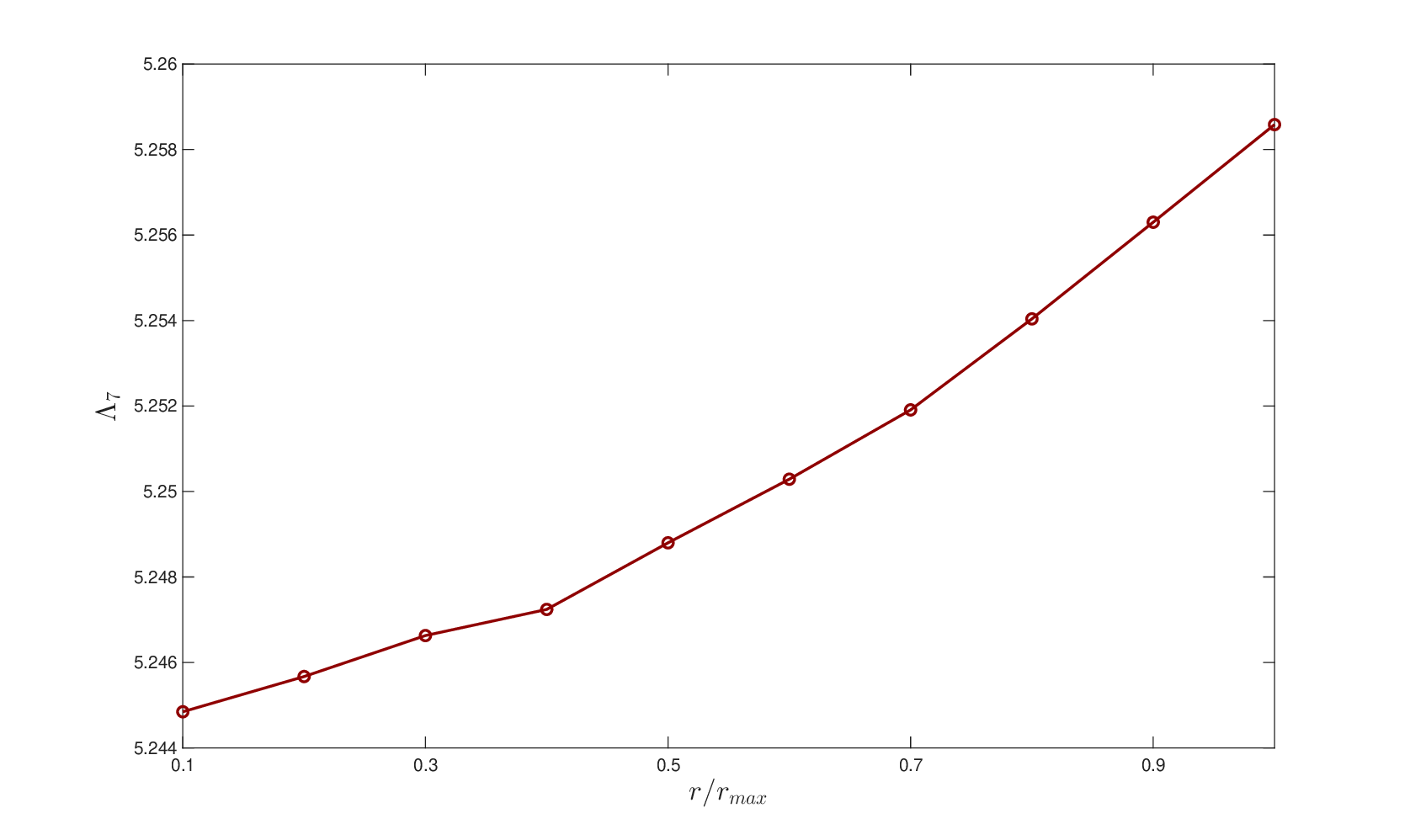}}
\caption{Lebesgue constants $ \Lambda_7 $ of discs centered at Chebyshev points. Even if the size of the discs varies, the Lebesgue constant remains rather stable (note the scale on the $y$-axis).} \label{fig:LebRadVary}
\end{figure}

Numerically, a further instance of stability 
emerges. In particular, the location of the centers of the discs appears to impact on the Lebesgue constant than their radius. To have an account of this, we set up the following experiment. 
Consider a unisolvent set of nodes $ \mathcal X = \{ \xi_1, \ldots, \xi_N \} $ and the discs $ K_1, \ldots, K_N $ (all sharing the same radius $ r $) centered at $ \mathcal X $. 
%
%
In Figure \ref{fig:LebRadVary} we compute $ \Lambda_7 $ as a function of $ r $, expressed as the ratio $ r / r_{\max} $, being $ r_{\max} $ the largest radius for which the discs are disjoint. Notice that, in this case, it is not guaranteed that the hypotheses of Proposition \ref{proposition:stability} are fulfilled. Nevertheless, we only observe a small change in the Lebesgue constant as the radius varies. 

This justifies, in the next numerical tests, the choice of a fixed radius $ r $ for all discs, easing the description of the supports $ B(\xi_i, r) $ by only the parameter $ \xi_i $ representing the center. Also, this radius will be taken in such a way that the hypotheses of Theorem \ref{thm:normLeb} are fulfilled, so that $ \Vert \Pi \Vert_{\mathrm{op}} = \Lambda_d $.


%

%
\subsubsection{Approximate Fekete supports and discrete Leja supports} \label{sect:FeketeLeja}

All the preceding techniques hinge on a collection of nodes for the construction of a unisolvent collection of compact sets. There is, however, another possibility, which also yields generally stable supports for interpolation. This procedure mimics Fekete and Leja constructions.

In analogy with the nodal case, we call \emph{Fekete supports} for $ \P_d (\Omega) $ the collection of compacts $ \{ K_1, \ldots, K_N \} $ that maximize the modulus of the determinant of the (normalized) Vandermonde matrix \eqref{eq:VdM}. While the exact problem is independent of the basis chosen for the polynomial space \cite{Bos91}, in the approximate case the selection of the supports depends on the basis considered for $ \P_d (\Omega) $, see  \cite[p. 17]{BDSV11}. Consistent with the related literature, we consider the basis detected in Section \ref{sect:choiceofthebasis} and, further, we normalize basis functions, as described in Remark \ref{rmk:normalizedbasis}. This is crucial when considering varying size supports \cite{BEFekete}.

Supports $ \{ K_1, \ldots, K_N \} $ that maximize the modulus of the determinant of $ \boldsymbol V $, indeed provide a controlled growth of the corresponding Lebesgue constant. 
In fact, since $ \Vert \widetilde{\ell}_{K_i} \Vert_\infty \leq 1 $ for $ i = 1, \ldots, N = \dim \P_d (\Omega) $, then
\begin{equation} \label{eq:estimateFekete}
\Lambda_d \leq \dim \P_d (\Omega) .
\end{equation}

\begin{figure}[!h]
   \centering
        {\includegraphics[width=5.9cm]{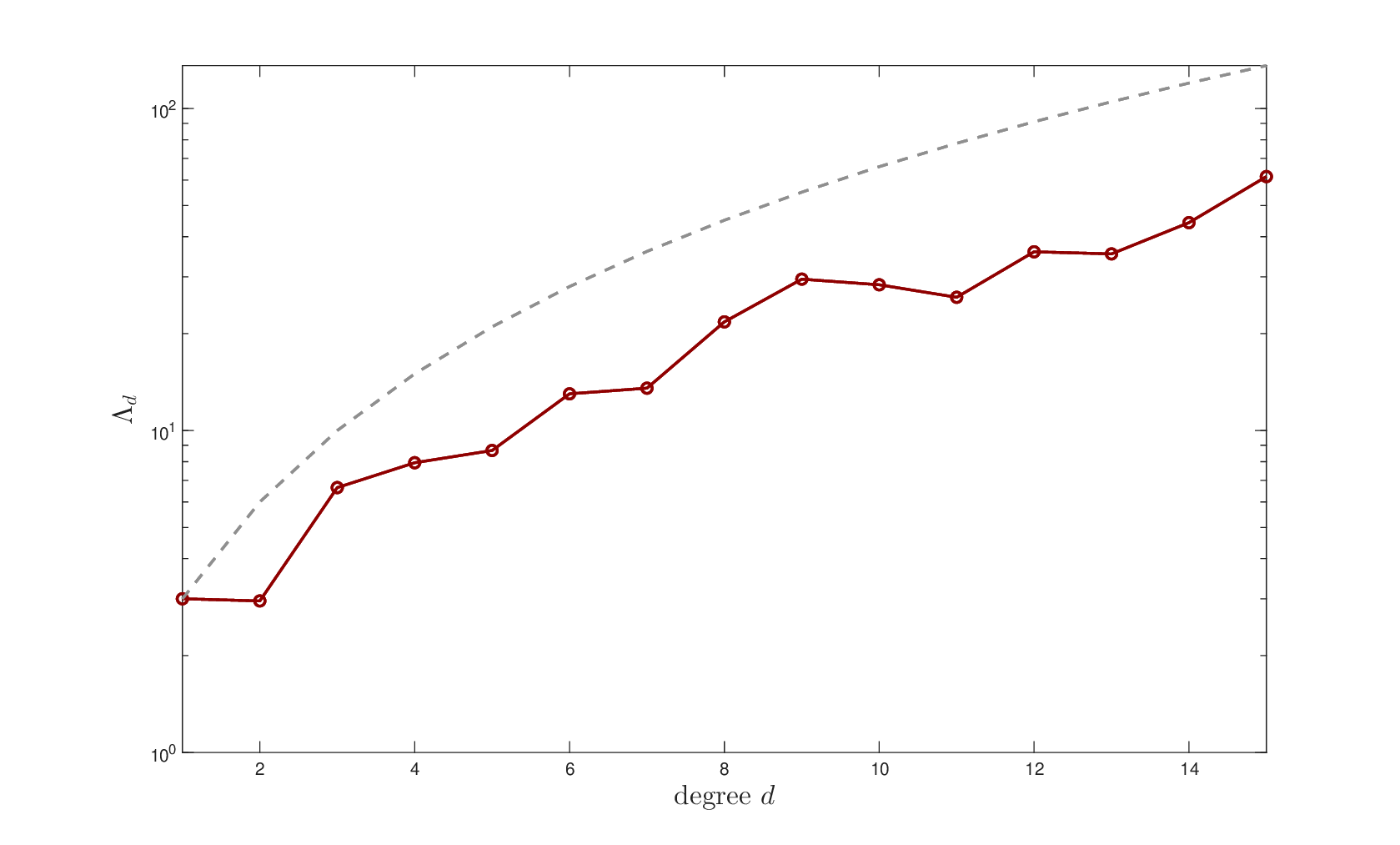}}
     \quad
        {\includegraphics[width=5.9cm]{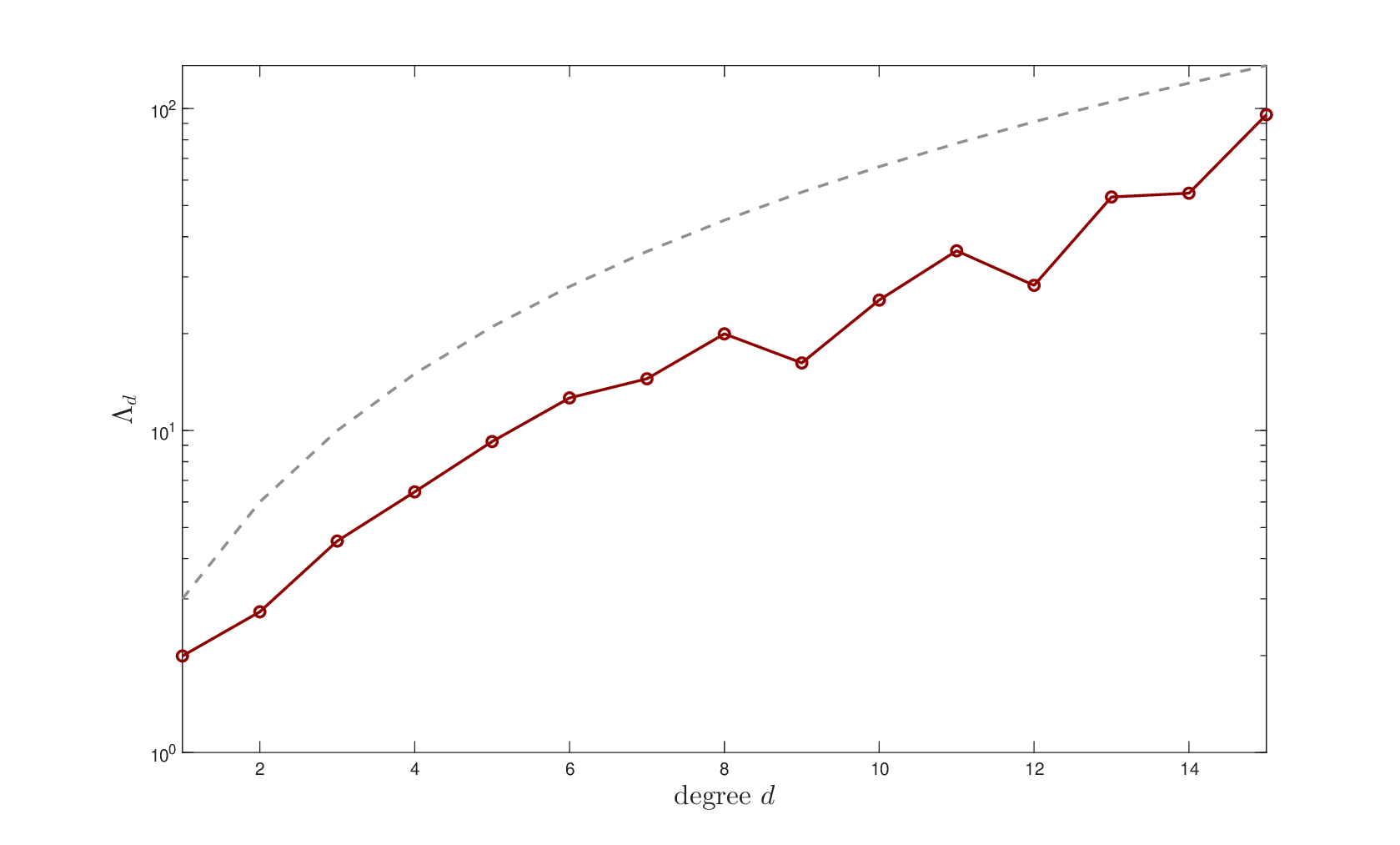}}     
\caption{Lebesgue constants for $ d = 1, \ldots, 15 $: approximate Fekete supports (left, solid red) and discrete Leja supports (right, solid red), compared with the theoretical bound (gray, dashed) given in \eqref{eq:estimateFekete}. Supports are extracted from a mesh containing about $ 8000 $ discs centered at a uniform grid.} \label{fig:LebAFS-DLS}
\end{figure}

The identification of Fekete supports is a hard problem \cite{BEFekete}, and algorithms for the approximate Fekete nodes have been designed by exploiting a greedy procedure on different bi and tri-dimensional domains \cite{BCLSV11, BSV12, SV09}. 

\begin{remark} 
From the theoretical point of view \cite[Sect. 3]{BDSV10}, the supports identified by such algorithms may violate the estimate \eqref{eq:estimateFekete}, as the search does not involve all possible $ K \subset \Omega $. On the numerical side, this effect has been observed in the simplicial framework \cite[Fig. 6]{PreprintBDEN}. 
\end{remark}

In the nodal case, similar to approximate Fekete points, are the {\it discrete Leja sequences} \cite{dMLeja}. Therefore, we call \emph{Leja (sequences of) supports} the arrays of supports $ \{ K_1, \ldots, K_N \} $, $ N = \dim \P_d (\Omega) $, obtained as follows. Once a basis is fixed, $ \{ b_1, \ldots, b_N \} $ for $ \P_d (\Omega) $, the first support $K_1$ of the sequence is
$$ K_1 = \arg \max_{K\subset \Omega} \frac{1}{\mu(K)}\left| \int_K b_1(x) \de x \right|. $$
The subsequent supports are chosen iteratively so that, at the $\ell$-th step of the procedure, the modulus of the determinant of the partial Vandermonde matrix 
$$ \left| \det \left[ \frac{1}{\mu(K_i)}\int_{K_i} b_j (x) \de x \right]_{i,j = 1}^\ell \right| $$ 
is maximized (the inductive definition ensures that $ \{ K_1, \ldots, K_{\ell-1} \} $ have been fixed). For greedy algorithms for the construction of nodal (discrete) Leja sequences see \cite{BDSV10,BDSV11}.

\begin{remark}  Leja sequences depend on (the ordering of) the basis chosen for the space of polynomials \cite[Section 3]{BDSV10}. However, once two bases $\{b_i\}_{i=1}^{\dim \bbP_d(\Omega)}$ and $\{g_i\}_{i=1}^{\dim \bbP_d(\Omega)}$ satisfy
$$ \spa \{ b_1,\dots,b_k\} = \spa \{ g_1,\dots,g_k\} \qquad  k=1,\dots,\dim \bbP_d(\Omega),$$
the corresponding (discrete) Leja sequences coincide \cite[p. 17]{BDSV11}. In fact, since such a condition holds for $ \mathcal{M}_d$ and $ \mathcal{T}_d$, the selected Discrete Leja sequences are the same for both bases.
\end{remark}

Those algorithms applied to supports provide both the {\it approximate Fekete supports (AFS)} and the discrete Leja supports (DLS). The greedy maximization of submatrix volumes by the QR factorization with column pivoting of the transposed Vandermonde matrix \eqref{eq:VdM} gives rise to the AFS, while the 
DLS corresponds to the greedy
maximization of nested square submatrix determinants (implemented by a standard LU factorization with row pivoting). The Lebesgue constants of the corresponding sets are shown in Figure \ref{fig:LebAFS-DLS} for the AFS (left) and for the DLS (right). The gray dashed line represents the bound \eqref{eq:estimateFekete}.

\subsection{Interpolation tests}

In Proposition \ref{prop:unisolvencepoints} and Theorem \ref{thm:algvar} we have provided two methodologies for constructing unisolvent sets for the interpolation problem \eqref{eq:areamatching}. The first hinges on points-based supports, whereas the second exploits a decomposition into algebraic varieties (or possibly orbits, as in Proposition \ref{thm:unisolvenceorbits}). Further, in Section \ref{sect:FeketeLeja}, we extended algorithms for the search of unisolvent and well-performing points to the case of diffused supports. We now compare these approaches in terms of the resulting interpolator.

Concerning points-based unisolvent sets, we shall consider either discs centered at quasi-random Halton points \cite{Halton} or at the \textquotedblleft Chebyshev\textquotedblright\ points of \cite{BX03} (for an optimized version, see \cite{MS19}), whose Lebesgue constant is depicted in Figure \ref{fig:errorbar} for the degree 7. Notice that this latter collection also fits the class of orbit-based supports. 

\begin{figure}[H]
\centering
\includegraphics[width=6.5cm]{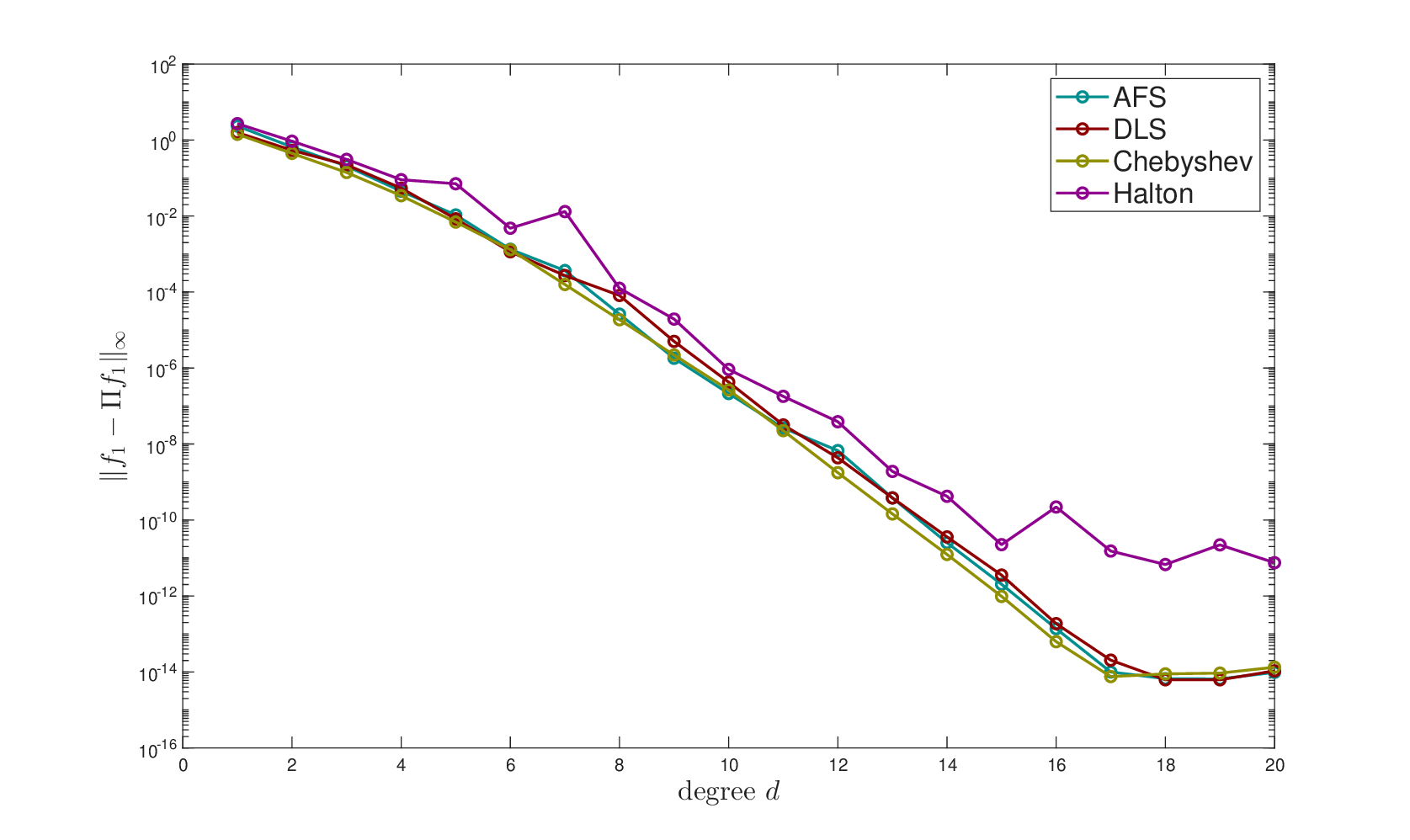}
\caption{Error $ \Vert f_1 - \Pi f_1 \Vert_\infty $ for $ d = 1, \ldots, 20 $, for different choices of supports. All the radii are constant and chosen to avoid overlaps.}
\label{fig:errf1}
\end{figure}

We consider the function
$$ f_1 (x,y) = e^x \sin(x+y).$$
%
In Figure \ref{fig:errf1} we depict the interpolation error $ \Vert f_1 - \Pi f_1 \Vert_\infty $, $ \Pi $ being defined as in \eqref{eq:Lagrangerepresentation}, with respect to the sup-norm, as a function of the total polynomial degree $ d $ of the interpolant $ \Pi $. All the strategies at play behave similarly, showing comparable convergence towards the exact solution. 

\begin{figure}[!h]
   \centering
        {\includegraphics[width=5.9cm]{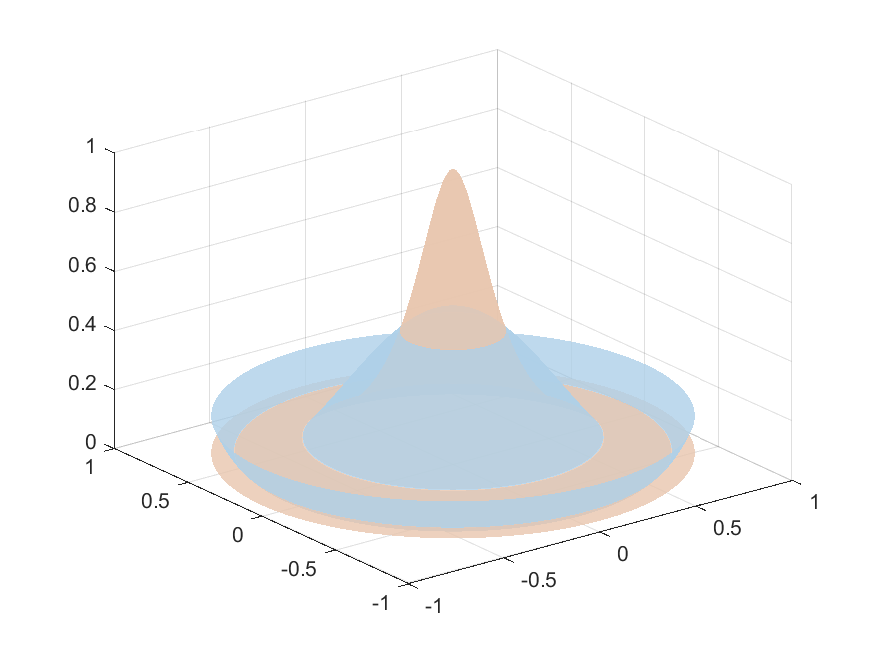}} \quad
        {\includegraphics[width=5.9cm]{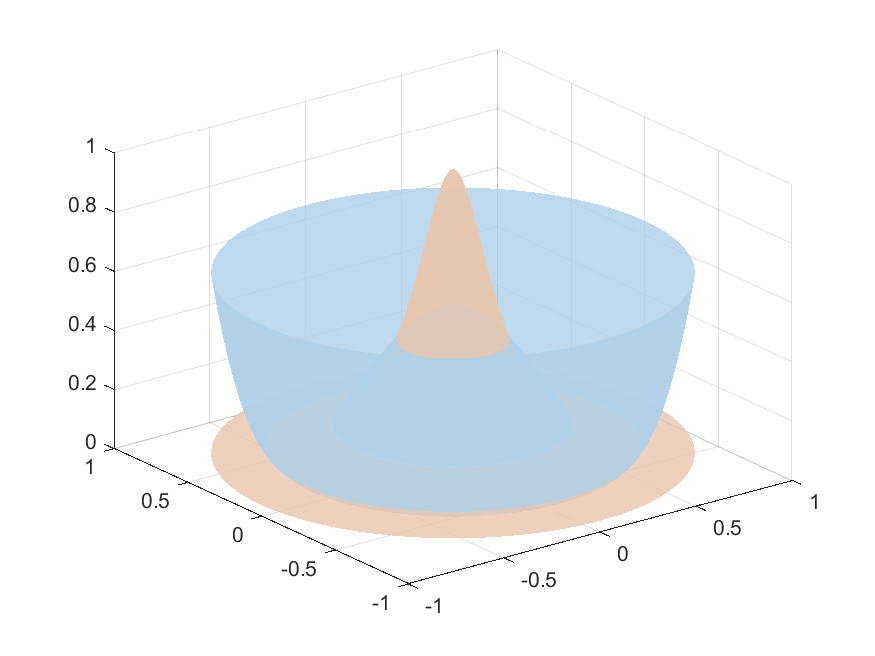}}     
\caption{The exact function $f_2$ and two interpolants. One based on discs with centers driven by Halton points (right), the other with radius based on a Chebyshev distribution (left), for $ d = 5 $.}
\label{fig:interpRunge}
\end{figure}

The situation becomes more interesting when considering a Runge-like function \cite{RothThesis}
$$ f_2 (x,y) = \frac{1}{25(x^2+y^2)+1} .$$
In this case, the choice of supports seriously affects the reliability of the interpolation. In accordance with the expectations, if supports are not carefully selected, the interpolant captures the function only away from the boundary, where instability becomes evident, as depicted in Figure \ref{fig:interpRunge}.

\begin{figure}[H]
\centering
\includegraphics[width=6.5cm]{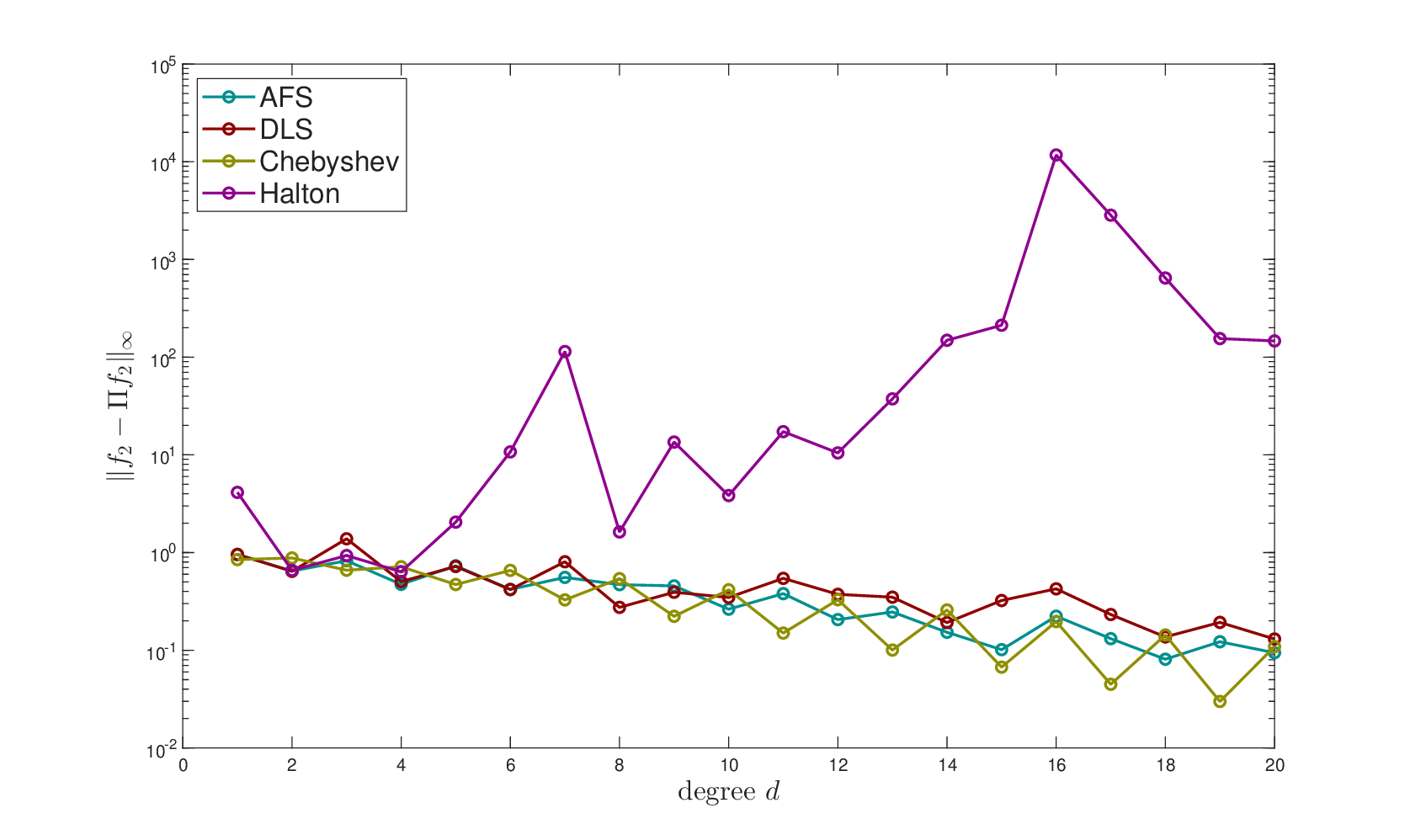}
\caption{Error $ \Vert f_2 - \Pi f_2 \Vert_\infty $ for $ d = 1, \ldots, 20 $, for different choices of supports. All the radii are constant and chosen to avoid overlaps.}
\label{fig:errf2}
\end{figure}

Figure \ref{fig:errf2} reports the errors $ \Vert f_2 - \Pi f_2 \Vert_\infty $ again with respect to the total degree $ d$. This indicates that the convergence towards the exact function is strongly dependent on the compact sets selected as supports of integration. In particular, it shows that the quantity $ \Lambda_d $ defined in Eq. \eqref{eq:Leb} fruitfully predicts the quality of the interpolation supports.

%

\section{Conclusions}

In this paper, we investigated the existence and features of a Lagrange-like polynomial interpolation scheme based on diffuse data. Such data is measured via integral quantities on compact supports. We showed that the properties of the corresponding projector $ \Pi $ closely resemble those of the nodal Lagrange interpolator. We introduced the concept of the Lebesgue constant $ \Lambda_d $, which relates the supports of integration to $ \Vert \Pi \Vert_{\mathrm{op}} $, where the uniform norm is considered. After determining a cost-effective method to compute $ \Lambda_d $, minimizing quadrature as much as possible, we employed three different strategies to construct unisolvent sets and assessed their efficacy. This also comes with a stability result that allows us to relate supports with some nodes, for which a large literature is available. As an application of this, we matched the estimated Lebesgue constants $ \Lambda_d $ with some interpolation errors coming from Runge-like problems. Numerical results show that the positioning of the compact supports, rather than their size, controls the features of the interpolator $ \Pi $, and that all the strategies proposed for the construction of unisolvent sets of compact supports may be used to extract productive results. In the future, we aim to consolidate the understanding of the trend of the Lebesgue constants by determining explicit bounds for their growth.

\section*{Acknowledgements}

This research has been accomplished within the thematic group on \lq\lq Teoria dell'Ap\-pros\-si\-ma\-zio\-ne e Applicazioni\rq\rq $\,$(TAA) of the Italian Mathematical Union and partially supported by GNCS-IN$\delta$AM. The third author have been funded by the Swiss National Science Foundation starting grant
``Multiresolution methods for unstructured data'' (TMSGI2\_211684).

\paragraph{Data availability}

All the numerical data shown in the present work have been produced by the authors. Open-source codes are available at the repository \\ \texttt{https://github.com/gelefant/InterpDISC} (for the interpolation test) and \\ \texttt{https://github.com/gelefant/DiscAFS\textunderscore DLS} (for the extraction of Fekete and Leja supports).

\printbibliography

\end{document}